\documentclass[reqno,a4paper,10pt]{amsart}
\usepackage{amsfonts,amsthm,amsmath,latexsym,bm,graphics,indentfirst,color}
\usepackage[hypertex]{hyperref}
\usepackage{geometry}
\newtheorem{theorem}{Theorem}[section]
\newtheorem{lemma}[theorem]{Lemma}
\newtheorem{proposition}[theorem]{Proposition}
\newtheorem{corollary}[theorem]{Corollary}
\newtheorem{remark}[theorem]{Remark}

\newcommand{\e}{\varepsilon}
\newcommand{\D}{\displaystyle}
\newcommand{\F}[2]{\frac{#1}{#2}}

\newcommand{\BE}{\begin{equation}}
\newcommand{\BEN}{\begin{equation*}}
\newcommand{\EE}{\end{equation}}
\newcommand{\EEN}{\end{equation*}}
\newcommand{\BL}{\begin{lemma}}
\newcommand{\EL}{\end{lemma}}
\newcommand{\BT}{\begin{theorem}}
\newcommand{\ET}{\end{theorem}}
\newcommand{\BP}{\begin{proposition}}
\newcommand{\EP}{\end{proposition}}
\newcommand{\BC}{\begin{corollary}}
\newcommand{\EC}{\end{corollary}}
\newcommand{\BR}{\begin{remark}}
\newcommand{\ER}{\end{remark}}

\allowdisplaybreaks

\begin{document}
\title[concentration solutions]{Concentration solutions to singularly  prescribed Gaussian and geodesic curvatures problem }

\author{Liping Wang}
\address{Liping Wang School of Mathematical Sciences, Shanghai Key Laboratory of Pure Mathematics and Mathematical Practice, East China Normal University, 200241, P.R. China}
\email{lpwang@math.ecnu.edu.cn}

\author{Chunyi Zhao}
\address{Chunyi Zhao School of Mathematical Sciences, Shanghai Key Laboratory of Pure Mathematics and Mathematical Practice, East China Normal University, 200241, P.R. China}
\email{cyzhao@math.ecnu.edu.cn}

\thanks{The research of the first author is supported  by NSFC 11671144. The research of Zhao is partially supported by NSFC 11971169.}

\keywords{ concentration, existence, exponential Neumann boundary condition }

\subjclass{}

\begin{abstract}
We consider the following Liouville-type equation with exponential Neumann boundary condition:
\[
\begin{cases}
  -\Delta\tilde u = \e^2 K(x) e^{2\tilde u}, & x\in D,  \\
  \D\F{\partial \tilde u}{\partial \bm n} + 1 = \e \kappa(x) e^{\tilde u},  \qquad &  x\in\partial D,
\end{cases}
\]
where $D\subset \mathbb R^2$ is the unit disc, $\e^2 K(x)$ and $\e \kappa(x)$ stand for the prescribed Gaussian curvature and the prescribed geodesic curvature of the boundary, respectively.    We prove the existence of concentration solutions if $\kappa(x) + \sqrt{K(x)+\kappa(x)^2}$ ($x\in\partial D$) has a  strictly local extremum point, which  is a total new result for exponential Neumann boundary problem.
\end{abstract}

\maketitle

\section{Introduction}

One of classical problems in conformal geometry is to prescribe the Gaussian curvature on a closed Riemannian surface $(\mathcal M,\tilde g)$. Denote by $g$ the conformal one, and by $e^{2v}$ the conformal factor. Then $v$ should comply with the elliptic equation that
\BEN
-\Delta_{\tilde g} v + \widetilde K(x) = K(x) e^{2v}.
\EEN
Here $\widetilde K$ and $K$ stands for the Gaussian curvature with respect to $\tilde g$ and $g$ respectively. After the pioneer works \cite{KW,KW1} by Kazdan and Warner, large amounts of literature are devoted to the solvability of this type of problems, see Chapter 6 in the book \cite{A}, where one can also find a comprehensive list of references.

In the case the boundary  $\partial \mathcal M \neq \emptyset$, other than the Gaussian curvature, it is natural to prescribe also the geodesic curvature on $\partial \mathcal M$. Let $\tilde \kappa$, $\kappa$ be the geodesic curvatures of the boundary with respect to $\tilde g$ and $g$ respectively. We are then led to the problem
\BE\label{12}
\begin{cases}
  -\Delta_{\tilde g}v + \widetilde K(x) =  K(x) e^{2v},\qquad & x\in \mathcal M,  \\
  \D\F{\partial v}{\partial \bm n} + \tilde \kappa(x) =  \kappa(x) e^{v},  \qquad &  x\in\partial \mathcal M.
\end{cases}
\EE
Brendle \cite{B} found a solution of  (\ref{12}) for constant $K$ and $\kappa$ by using the method of flows. When $K$ and $\kappa$ are not constants, Cherrier \cite{C} proved the existence of a solution for not big curvatures.
Recently, L\'{o}pez-Soriano, Malchiodi and Ruiz \cite{LMR} considered the negative Gaussian curvature case, i.e $K<0$,  and they derived some existence results using a variational approach.

The higher-dimensional analogue of this question is to prescribe scalar curvature of a manifold and mean curvature of the boundary.
Moreover, the case of zero scalar curvature has been widely studied, see for instance \cite{ALM, CXY, DMO, E, HL1, WZ} and the references therein.

When $\mathcal M$ is the standard Euclidean disk $D \subset \mathbb R^2$, this problem may be recognised as a kind of Nirenberg problem to surfaces with boundary.  The case $\kappa = 0$ has been treated by Chang-Yang \cite{CY}, which proved the existence of a solution to (2) under the assumption that $K$ is positive somewhere. Moreover,
the case $K = 0$ was first treated in Chang-Liu \cite{CL} under a suitable conditions on critical points of the geodesic curvature.   Also in \cite{LH}  a solution is found under some symmetric condition. Additionally, a blow up analysis has been given in \cite{GJ}. \par

For the nonconstant curvatures case,  the first work \cite{CR} address the existence by setting the problem in a variational framework also under some symmetric assumptions.
Recently, Jevnikar and etc \cite{JLMR} gave a complete blow up analysis of  the problem (\ref{12}), where $K$ and $\kappa$ are not constant. Precisely, they considered the problem
\BEN
\begin{cases}
  -\Delta u_n = 2 K_n e^{u_n}, \qquad &\text{in }D,\\
  \F{\partial u_n}{\partial\nu} +2 =2\kappa_n e^\F{u_n}{2}, &\text{on }\partial D,
\end{cases}
\EEN
where $K_n \to K$ in $C^2(\overline { D})$, $\kappa_n \to\kappa$ in $C^2(\partial{ D})$ and $K,\kappa$ are not all zeros. If $u_n$ blows up with bounded mass, i.e.
\BEN
\int_{\mathbb D} e^{u_n} + \int_{\partial \mathbb D} e^\F{u_n}{2} \leq C <+\infty,
\EEN
then  the blow-up must occur at a critical point of $\varphi(\xi) = H(\xi) + \sqrt{H(\xi)^2 + K(\xi)}$, where $H$ is the harmonic extension of $\kappa$. When finishing our writing the other day, we find that Battaglia-Medina-Pistoia \cite{BMP} considered the inverse problem to \cite{JLMR}. They constructed a family of solutions with convergent curvatures $K_\e\to K$ and $\kappa_\e\to\kappa$ as $\e$ goes to $0$, which blows up at one boundary critical point of $\varphi(\xi)$ under some non-degeneracy conditions. Note that in their settings, $K$ is not identical to $0$.

In this paper, we consider the following perturbation problem
\BE\label{1}
\begin{cases}
  -\Delta\tilde u = \e^2 K(x) e^{2\tilde u}, & x\in D,  \\
  \D\F{\partial \tilde u}{\partial \bm n} + 1 = \e \kappa(x) e^{\tilde u},  \qquad &  x\in\partial D,
\end{cases}
\EE
where $D\subset \mathbb R^2$ is the unit disc (\textbf{centered at $(0,1)$ for convience}), $K(x)$ and $\kappa(x)$ are all positive smooth functions on $\overline D$ and $\partial D$ respectively. Here $\e>0$ is a small parameter. Note that $\e^2K(x)$ and $\e\kappa(x)$ stand for the prescribed Gaussian curvature and the prescribed geodesic curvature of the boundary.  We mention that the problem (\ref{1}) is the perturbation case of zero curvatures, which do not be included in \cite{BMP} and \cite{JLMR}.  Note that the universal constant $ 1$ is the original geodesic curvature of the boundary $\partial  D$. Our main result are stated as follows.

\BT\label{T1}
Assume that the function $\kappa(x) + \sqrt{K(x)+\kappa(x)^2}$ defined on $\partial D$ admits a local extremum point. Then there is an $\e_0>0$ such that for any small $0<\e<\e_0$, the problem (\ref{1})
has a family of single boundary bubbling solutions $u_\e$. Moreover,
\BEN
\int_D \e^2 K(x) e^{2\tilde u_\e} + \int_{\partial D}\e \kappa(x) e^{\tilde u_\e} = 2\pi.
\EEN
\ET

\BR
As long as $\kappa(x) + \sqrt{K(x)+\kappa(x)^2}$ is not a constant on $\partial D$,  the problem (\ref{1}) admits single boundary bubbling solutions.
\ER

The proof of Theorem \ref{T1} is based on the Lypunov-Schmidt reduction as performed in \cite{DPM}. We first make use of the classification result on $\mathbb R^2_+$ \cite{LZ} to build up the approximation solution. For the non-degeneracy of standard bubble, the disk $D$ is transformed to a sphere cap by stereographic projection, and then  the non-degeneracy follows from the result in \cite{HL}.
For the solvability of corresponding linearized problem,  there seems a gap in getting (53) from (75) in \cite{DPM}. In this paper we give a new  proof instead.\par

This paper is organized as follows. In Section \ref{s2}, an ansatz of the solution is given. Section \ref{s3} is devoted to the invertibility of linearized operator. In Section \ref{s4}, the nonlinear problem is solved. Variational reduction is then showed in Section \ref{s5}. In Section \ref{s6} the theorem \ref{T1} is proved. Some computation is listed in the appendix section \ref{s7}.

\section{Ansatz}\label{s2}
In this section, we will introduce the approximation solution and give the ansatz of the solution.\par
Let $u(y)=\tilde u(\e y) + 2\ln\e$ and $D_\e = D/\e$. Then the equation (\ref{1}) is equivalent to
\BE\label{2}
\begin{cases}
  -\Delta u =  K(\e y) e^{2u}, & y\in D_\e,  \\
  \D\F{\partial u}{\partial \bm n} + \e  = \kappa(\e y) e^{u},  \qquad &  y\in\partial D_\e.
\end{cases}
\EE
For $\e \to 0$, it is easy to understand that the limit problem near a boundary point $\xi=\e \xi'$ is formally like
\BEN
\begin{cases}
  -\Delta u =  K(\xi) e^{2u}, & y\in \mathbb R^2_+,  \\
  \D\F{\partial u}{\partial \bm n}  = \kappa(\xi) e^{u},  \qquad &  y\in\partial \mathbb R^2_+.
\end{cases}
\EEN
\par

To introduce the approximation solution, we recall that the unique solution (\cite{LZ}) to the following half plane problem
\BEN
\begin{cases}
  -\Delta u = a e^{2u}, \quad  & x\in \mathbb R_+^2, \\
  \D \F{\partial u}{\partial \bm n} = b e^u, & x\in\partial\mathbb R_+^2,
\end{cases}\qquad
(a,b>0 \text{ are two constants})
\EEN
is
\BEN
\widetilde U_0(x) = \ln \F{2\lambda}{\sqrt{a}\left(\lambda^2+(x_1-s)^2+(x_2+\F{b}{\sqrt a}\lambda)^2\right)}, \quad (s \text{ and }\lambda\text{ are two arbitrary parameters}).
\EEN
From the above classification result,   the first approximation we select is  that
\BEN
U_0(x) = \ln \F{2\lambda}{\sqrt{K(\xi)}\left(\lambda^2\e^2+\left|x-\xi-\F{\kappa(\xi)}{\sqrt{K(\xi)}}\lambda\e\bm n(\xi)\right|^2\right)}.
\EEN
Here $\xi\in\partial D$ is undetermined and $\bm n(\xi)$ is the unit out normal at the point $\xi$, while $\lambda$ may be any universal constant. Obviously $U_0$ satisfies the problem
\BEN
\begin{cases}
  -\Delta U_0 = \e^2 K(\xi) e^{2U_0} , \quad & x\in \mathbb R_+^2, \\
  \D \F{\partial U_0}{\partial \bm n} = \e \kappa(\xi) e^{U_0},& x\in\partial\mathbb R_+^2.
\end{cases}
\EEN
It is found that the function $U_0$ approximate the equation (\ref{1}) near $\xi$, but far from the boundary  condition.
As usual, we will modify $U_0$ furthermore.
On the boundary, one may find that
\begin{align}
&\ \e \kappa(\xi)e^{U_0} - \F{\partial U_0}{\partial\bm n} \nonumber \\
=&\ \F{2\kappa(\xi)\lambda \e } {\sqrt{K(\xi)}\left(\lambda^2\e^2+\left|x-\xi-\F{\kappa(\xi)}{\sqrt{K(\xi)}}\lambda\e\bm n(\xi)\right|^2\right)}
 + \F{2\left(x-\xi-\F{\kappa(\xi)}{\sqrt{K(\xi)}}\lambda\e\bm n(\xi)\right)\cdot \bm n(x)}{\lambda^2\e^2+\left|x-\xi-\F{\kappa(\xi)}{\sqrt{K(\xi)}}\lambda\e\bm n(\xi)\right|^2} \nonumber\\
=&\ \F{2(x-\xi)\cdot\bm n(x)}{\lambda^2\e^2+\left|x-\xi-\F{\kappa(\xi)}{\sqrt{K(\xi)}}\lambda\e\bm n(\xi)\right|^2} + \F{2\kappa(\xi)\lambda\e(1-\bm n(\xi)\cdot\bm n(x))}{\sqrt{K(\xi)}\left(\lambda^2\e^2+\left|x-\xi-\F{\kappa(\xi)}{\sqrt{K(\xi)}}\lambda\e\bm n(\xi)\right|^2\right)} \nonumber\\
=&\ \F{2(x-\xi)\cdot\bm n(x)}{|x-\xi|^2} + \left(\F{2(x-\xi)\cdot\bm n(x)}{\lambda^2\e^2+\left|x-\xi-\F{\kappa(\xi)}{\sqrt{K(\xi)}}\lambda\e\bm n(\xi)\right|^2}-\F{2(x-\xi)\cdot\bm n(x)}{|x-\xi|^2}\right) \nonumber\\
&\qquad + \F{2\kappa(\xi)\lambda\e(1-\bm n(\xi)\cdot\bm n(x))}{\sqrt{K(\xi)}\left(\lambda^2\e^2+\left|x-\xi-\F{\kappa(\xi)}{\sqrt{K(\xi)}}\lambda\e\bm n(\xi)\right|^2\right)} \nonumber\\
=&\ \F{2(x-\xi)\cdot\bm n(x)}{|x-\xi|^2} + I_1(x) +I_2(x). \label{10}
\end{align}
Set $\mathfrak D(\xi) = \F{\kappa(\xi)}{\sqrt{K(\xi)}}$ in what follows for convenience. It is checked that
\begin{align*}
  |I_1(x)| =&\ \left|2\lambda \e \F{(x-\xi)\cdot\bm n(x)\left[2\mathfrak{D}(\xi)(x-\xi)\cdot\bm n(\xi) - (1+\mathfrak{D}(\xi)^2)\lambda \e \right]}{|x-\xi|^2\left(\lambda^2\e^2+\left|x-\xi- \mathfrak{D}(\xi)\lambda\e\bm n(\xi)\right|^2\right)}\right| \\
  \leq &\ C\lambda\e + C \lambda^2\e^2 \F{1}{\lambda^2\e^2+\left|x-\xi- \mathfrak{D}(\xi)\lambda\e\bm n(\xi)\right|^2},
\end{align*}
since $\F{2(x-\xi)\cdot\bm n(x)}{|x-\xi|^2}\equiv 1$ for $x\in\partial D$.
So it is easily to get, by assuming $\xi=(0,0)$ from the symmetry and using $x-\xi = \lambda\e z$ , that
\BEN
\left|\int_{\partial D} I_1\right| \leq C\lambda \e\int_{\partial D_{\lambda\e}} \F{\mathrm d z}{1+|z -\mathfrak{D}(0)\bm n(0)|^2}  + C\lambda \e \leq C\lambda\e.
\EEN
Also it holds that
\BEN
I_2(x) = O(\lambda\e), \qquad \left|\int_{\partial D} I_2\right| \leq C\lambda\e.
\EEN
Next, we define a constant $d=O(1)$ by
\BEN
  d \int_{\partial D} \F{\lambda^2\e^2 \mathrm d x}{\lambda^2\e^2 + |x-\xi|^2}= \int_{\partial D} (I_1 + I_2) .
\EEN
Thus it is easy to see, owing to $\int_{\partial D} \left(\F{2(x-\xi)\cdot\bm n(x)}{|x-\xi|^2} - 1\right) \mathrm d x = 0$ , that
\BEN
\int_{\partial D} \left( \e \kappa(\xi)e^{U_0} - \F{\partial U_0}{\partial\bm n} - 1 - \F{\lambda^2\e^2 d}{\lambda^2\e^2 + |x-\xi|^2} \right) \mathrm dx  = 0.
\EEN
So one may define a function $H_0(x)$ satisfying
\BEN
\begin{cases}
-\Delta H_0 = 0,  &  x\in D , \\
\D \F{\partial H_0}{\partial\bm n} =\e \kappa(\xi)e^{U_0} - \F{\partial U_0}{\partial\bm n} - 1 - \F{\lambda^2\e^2 d}{\lambda^2\e^2 + |x-\xi|^2}, \quad & x\in\partial D.
\end{cases}
\EEN
Of course $H_0$ is unique up to a constant.
Obviously (\ref{10}) implies that
\BEN
\F{\partial H_0}{\partial\bm n} =  I_1 + I_2 - \F{\lambda^2\e^2 d}{\lambda^2\e^2 + |x-\xi|^2},  \qquad  x\in\partial D.
\EEN
And it holds that, for any $p>1$,
\begin{align*}
\|I_1\|_{L^p(\partial D)} &\leq C\lambda\e + C\lambda^2\e^2 \left(\int_{\partial D} \F{\mathrm dx}{\left(\lambda^2\e^2+\left|x-\xi- \mathfrak{D}(\xi)\lambda\e\bm n(\xi)\right|^2\right)^p }\right)^\F{1}{p}=O(\e^\F{1}{p}), \\
\|I_2\|_{L^p(\partial D)} &= O(\e^\F{1}{p}), \qquad \left\|\F{\lambda^2\e^2 d}{\lambda^2\e^2 + |x-\xi|^2}\right\|_{L^p(\partial D)} =O(\e^\F{1}{p}).
\end{align*}
Thus by $L^p$ theory we have, for $0<s<\F{1}{p}$, that
\BEN
\|\nabla H_0\|_{W^{s,p}(B)} = O(\e^\F{1}{p}).
\EEN
In what follows, we always choose $H_0$ such that $\int_D H_0 = 0$, which means that
$$H_0(x) =  O(\e^\alpha), \qquad (\alpha=1/p \in (0,1)\text{ is arbitrary})$$

We then choose the second approximation $U(x) = U_0(x) + H_0(x)$. It is easy to see that
\BE\label{3}
\begin{cases}
-\Delta U = \e^2 K(\xi) e^{2U_0},  &  x\in D , \\
\D\F{\partial U}{\partial\bm n} + 1 +  \F{d\lambda^2\e^2}{\lambda^2\e^2 + |x-\xi|^2}  =\e \kappa(\xi)e^{U_0}, \qquad & x\in\partial D.
\end{cases}
\EE
We will seek a solution to (\ref{2}) of the form $V(y)+\phi(y)$ where $V(y) = U(\e y) + 2\ln \e$. Thus the problem can be stated as to find a solution $\phi$ of
\BEN
\begin{cases}
-\Delta \phi - 2K(\e y)e^{2V}\phi =R_1(y) + K(\e y) e^{2V}(e^{2\phi}-1-2\phi),\quad  & \text{in } \ D_{\e} , \\
\D \F{\partial\phi}{\partial\bm n} - \kappa(\e y) e^V \phi =R_2(y) + \kappa(\e y) e^{V}(e^{\phi}-1-\phi), & \text{on }\partial D_{\e},
\end{cases}
\EEN
where the error terms
\BEN
 R_1(y)=\Delta V + K(\e y)e^{2V} , \qquad R_2(y)=- \F{\partial V}{\partial\bm n} - \e  + \kappa(\e y)e^V.
\EEN

\BL\label{lem1}
Let $\delta>0$ be small and fixed. In $B_\F{\delta}{\e}(\xi')\cap \partial D_\e$, we have
\begin{align*}
 R_1(y) &= \F{4\lambda^2}{(\lambda^2+|y-\xi'-\mathfrak D(\xi)\lambda\bm n(\xi)|^2)^2} \left[O(\e |y-\xi'|)+O(\e^\alpha)\right], \\
 R_2(y) &= \F{\mathfrak D(\xi) 2\lambda}{(\lambda^2+|y-\xi'-\mathfrak D(\xi)\lambda\bm n(\xi)|^2)} \left[O(\e |y-\xi'|)+O(\e^\alpha)\right] + \e \F{d\lambda^2}{\lambda^2 + |y-\xi'|^2},
\end{align*}
while in $B_\F{\delta}{\e}(\xi')^c\cap \partial D_\e$,
\BEN
R_1(y) = O(\e^4), \qquad R_2(y) = O(\e^2).
\EEN
\EL

\begin{proof}
Direct computation shows that, for $|y-\xi'|\leq \F{\delta}{\e}$,
\begin{align*}
  R_1(y)=&\  \Delta V + K(\e y)e^{2V}  = K(\e y) e^{2(U_0(\e y) + H_0(\e y) + 2\ln\e)} - K(\xi) e^{2(U_0(\e y) + 2\ln\e)}   \\
  =&\ \F{4\lambda^2}{K(\xi)(\lambda^2+|y-\xi'-\mathfrak D(\xi)\lambda\bm n(\xi)|^2)^2} (K(\e y)e^{2H_0(\e y)} - K(\xi)) \\
  =&\ \F{4\lambda^2}{(\lambda^2+|y-\xi'-\mathfrak D(\xi)\lambda\bm n(\xi)|^2)^2} \left(\F{K(\e y)}{K(\xi)}e^{O(\e^\alpha)} - 1\right) \\
  =&\ \F{4\lambda^2}{(\lambda^2+|y-\xi'-\mathfrak D(\xi)\lambda\bm n(\xi)|^2)^2} \left[O(\e |y-\xi'|)+O(\e^\alpha)\right],
\end{align*}
and
\begin{align*}
  R_2(y)=&\  - \F{\partial V}{\partial\bm n} - \e + \kappa(\e y)e^V = \kappa(\e y)e^V - \kappa(\xi)e^{U_0(\e y)+2\ln\e} + \e \F{d\lambda^2}{\lambda^2 + |y-\xi'|^2} \\
  =&\ \F{\mathfrak D(\xi) 2\lambda}{(\lambda^2+|y-\xi'-\mathfrak D(\xi)\lambda\bm n(\xi)|^2)} \left[O(\e |y-\xi'|)+O(\e^\alpha)\right] + \e \F{d\lambda^2}{\lambda^2 + |y-\xi'|^2}.
\end{align*}
\par
Other estimates of the lemma are easily checked.
\end{proof}

In the last of this section, we build the non-degeneracy of the standard bubble $\widetilde U_0$, which plays an important role thereinafter. For any positive numbers $a$, $b$, $\lambda$, let
\begin{equation}\label{basis}
z_0(x)=\frac1{2\lambda}-\D\F{\lambda+\frac{b}{\sqrt{a}}(x_2 + \frac{b}{\sqrt{a}}\lambda)}{\lambda^2+x_1^2+(x_2 + \frac{b}{\sqrt{a}}\lambda)^2}, \qquad z_1(x)=\D\F{x_1}{\lambda^2+x_1^2+(x_2 + \frac{b}{\sqrt{a}}\lambda)^2}, \qquad x=(x_1, x_2)\in\mathbb R^2_+.
\end{equation}
Then the following non-degeneracy holds.

\begin{lemma}\label{sl}
Any bounded solutions of
\begin{equation}\label{linear}
\begin{cases}
\Delta \phi  +\D\F{8\lambda}{\left(\lambda^2+x_1^2+(x_2 + \frac{b}{\sqrt{a}}\lambda)^2\right)^2}\phi=0,  \qquad  &\text{in }  \mathbb{R}^2_+,\\
\D \F{\partial\phi}{\partial\bm n} - \D\F{2b\lambda}{\sqrt{a}\left(\lambda^2+x_1^2+( \frac{b}{\sqrt{a}}\lambda)^2\right)}\phi =0,  \qquad &\text{on }  \partial\mathbb{R}^2_+
\end{cases}
\end{equation}
is a linear combination of $z_0$ and $z_1$.
\end{lemma}
\begin{proof}
let $\Pi$ be the stereographic projection from the unit sphere in $\mathbb{R}^3$ centered at $(0, -\frac{b}{\sqrt{a}}\lambda, 0)$ onto $\mathbb{R}^2$. More specifically, let $\eta=(\eta_1,\eta_2,\eta_3)$ be the coordinates of $\mathbb{R}^3$ taking $(0, -\frac{b}{\sqrt{a}}\lambda, 0)$ as its origin and  $x=(x_1,x_2)$ be the coordinates of $\mathbb{R}^2$, we have
\[
\eta_1=\frac{2x_1}{\lambda^2+x_1^2+(x_2 + \frac{b}{\sqrt{a}}\lambda)^2}, \qquad \eta_2=\frac{2(x_2 + \frac{b}{\sqrt{a}}\lambda)}{\lambda^2+x_1^2+(x_2 + \frac{b}{\sqrt{a}}\lambda)^2},\qquad \eta_3=\frac{x_1^2+(x_2 + \frac{b}{\sqrt{a}}\lambda)^2-\lambda^2}{\lambda^2+x_1^2+(x_2 + \frac{b}{\sqrt{a}}\lambda)^2}.
\]
Let $\Sigma=\Pi^{-1}(\mathbb{R}^2_+)$. It is a spherical cap on $\mathbb S^2$.

Assume $\phi$ is a bounded solution of problem (\ref{linear}), we define a function $\Phi(\eta)$ on $\Sigma$ by
\[
\phi(x)=\Phi(\eta)\frac{2\lambda}{\lambda^2+x_1^2+(x_2 + \frac{b}{\sqrt{a}}\lambda)^2},
\]
then
\begin{equation*}
\begin{cases}
\Delta \Phi  + 2\Phi=0 , \\
\D \F{\partial\Phi}{\partial \nu} - \D\F{b}{\sqrt{a}} \Phi =0.
\end{cases}
\end{equation*}
Using Proposition 3.2 in \cite{HL}, we get the desired result.
\end{proof}

\section{The Linearized Operator}\label{s3}

In the section, the invertibility of the linearized operator is studied. The main result is the solvability of the following linear problem. For given $f$ and $h$,  find $(\phi, c_1)$ such that
\BEN
\begin{cases}
-\Delta \phi - 2K(\e y)e^{2V}\phi = f+ c_1\chi_\e Z_{1,\e},  \qquad & \text{in } D_{\e} , \\
\D \F{\partial\phi}{\partial\bm n} - \kappa(\e y) e^V \phi = h , \quad & \text{on } \partial D_{\e} \medskip\\
\D\int_{D_{\e}} \chi_\e Z_{1,\e}\phi=0.
\end{cases}
\EEN
where $f\in L^\infty(D_{\e})$, $h \in L^\infty(\partial D_{\e})$ and $Z_{0,\e}$, $Z_{1,\e}$, $\chi_\e$ are defined as follows. From now on we denote  $z_0$, $z_1$ with $a=K(\xi)$, $b=\kappa(\xi)$ in (\ref{basis}). Around the point $\xi' = \xi/\e \in \partial D_{\e}$, we consider a smooth change of variables
\[
F_\e(y)=\frac1\e F(\e y),
\]
where $F$: $B_\rho(\xi)\rightarrow M$ is a diffeomorphism and $M$ is an open neighborhood of the origin such that $F(B\cap  B_\rho(\xi))=\mathbb{R}^2_+\cap M$, $F(\partial B\cap  B_\rho(\xi))=\partial\mathbb{R}^2_+\cap M$. Also we can choose $F$ preserving area. Define
\[
Z_{0,\e}(y)=z_0(F_\e(y)), \qquad \qquad Z_{1,\e}(y)=z_1(F_\e(y)).
\]
Next, we select a large but fixed number $R_0$ and non-negative smooth function $\chi: \mathbb{R} \rightarrow \mathbb{R}$ such that $\chi(r)=1$ for $r \le R_0$ and $\chi(r)=0$ for $r \ge R_0+1, 0\le \chi\le 1$. Then set
\[
\chi_\e (y)=\chi(F_\e (y)).
\]

For simplicity, let us denote $W_1(y)=2K(\e y)e^{2V}$ and $W_2(y)= \kappa(\e y) e^V$, then
\begin{equation*}
\begin{array}{ll}
& W_1(y)=\D \frac{2K(\e y)}{K(\xi)}\D\F{4\lambda^2}{\left(\lambda^2+|y-\xi'-\mathfrak D(\xi)\lambda\bm n(\xi)|^2\right)^2}\left[1+O(\e|y-\xi'|)+O(\e^\alpha)\right],\\
& W_2(y)=\D \frac{\kappa(\e y)}{\sqrt{K(\xi)}}\D\F{2\lambda}{\lambda^2+|y-\xi'-\mathfrak D(\xi)\lambda\bm n(\xi)|^2}\left[1+O(\e|y-\xi'|)+O(\e^\alpha)\right].
\end{array}
\end{equation*}
for any $\alpha \in (0, 1)$.

For $f\in L^\infty(D_{\e})$ and $h \in L^\infty(\partial D_{\e})$ we define two norms
\begin{align*}
\|f\|_{**, D_{\e}}&=\sup_{y\in D_{\e}}
 |f(y)|\cdot \left(1+|y-\xi'|^{2}\ln^3 (1+ |y-\xi'|)\right), \\
 \|h\|_{*,\partial D_{\e}}&=\sup_{y\in \partial D_{\e}}
 |h(y)|\cdot \left(1+|y-\xi'|\ln^3 (1+ |y-\xi'|)\right).
\end{align*}

\begin{lemma}\label{barrier}
Choose $R_1$ large enough such that for  any  $\e>0$ small enough, there exists a smooth and positive function
\[
\psi: D_{\e}\setminus B_{R_1}({\xi}')\rightarrow \mathbb{R}
\]
so that
\begin{eqnarray*}
&\D -\Delta \psi -W_1\psi \ge  \frac{1}{1+|y-\xi'|^{2}\ln^3 (1+ |y-\xi'|)},  \qquad &\text{in} \quad  D_{{\e}}\setminus B_{R_1}({\xi}'),\\
&\D \F{\partial\psi}{\partial\bm n}
 -W_2 \psi \ge  \frac{1}{1+|y-\xi'|\ln^3 (1+ |y-\xi'|)},  \qquad &\text{on} \quad  \partial D_{{\e}}\setminus B_{R_1}({\xi}'),\\
& \psi>0,   &\text{in} \quad  D_{{\e}}\setminus B_{R_1}({\xi}'),\\
& \psi \ge 1,  &  \text{on} \quad   D_{{\e}}\cap \partial B_{R_1}({\xi}').
\end{eqnarray*}
The positive constants  $C$ is independent of $\e, R_1$ and $\psi$ is bounded uniformly
\[
0< \psi \le C \qquad \text{in} \quad  D_{{\e}}\setminus B_{R_1}({\xi}').
\]
\end{lemma}

\begin{proof}
Without loss of generality we may assume that $\xi$ is the original point.  We take
\[
\psi_1=\D\F{-y_2}{r\ln^3 r}, \qquad y=(y_1,y_2), \qquad  r=|y-\mathfrak D(\xi)\lambda\bm n(\xi)|,
\]
then
\[
\nabla \psi_1 =\frac{(\ln^3 r + 3\ln^2 r)y_2}{r^2\ln ^6 r}\left(\frac{y_1}r, \frac{y_2 + \lambda \mathfrak D(\xi)}r    \right) - \left(0,  \frac1{r\ln^3r}\right).
\]
Direct computation gives that
\begin{eqnarray*}
&-\Delta \psi_1 -W_1\psi_1 =O\left(\D\F{1}{r^2 \ln ^3 r}\right)\\
& \D \F{\partial\psi_1}{\partial\bm n}=\nabla \psi_1 \cdot (\e y_1, \e y_2 -1)= \D\F{\e y_2 y_1^2 + (\e y_2 - 1)y_2^2}{r^3\ln^3 r}  + \F{1-\e y_2}{r\ln^3 r} + O\left(\frac1{r\ln^4 r}  \right)
\end{eqnarray*}
for $R_1$ large  and   on the boundary
\[
W_2\psi_1=O\left(\frac1{r^2\ln^3 r}\right).
\]
Let
\[
\psi_2=1-\frac1{\ln r},
\]
then
\[
\nabla \psi_2 = \frac1{r\ln^2 r}\left(\frac{y_1}r,  \frac{y_2 + \lambda \mathfrak D(\xi)}r    \right)
\]
from which we can obtain
\[
-\Delta \psi_2 -W_1\psi_2=\frac2{r^{2}\ln^3 r } + O\left(\frac1{r^{4}}\right),
\]
and
\[
\D \F{\partial\psi_2}{\partial\bm n}-W_2\psi_2=\frac{y_2}{r^{2}\ln^2 r} +O(\frac1{r^2}).
\]
Now we can define $\psi_0 = \psi_1 + C\psi_2$  then
\[
-\Delta \psi_0 -W_1\psi =\frac{2C}{r^{2}\ln^3 r } + O\left(\D\F{1}{r^2 \ln ^3 r}\right) \ge \frac{1}{r^{2}\ln^3 r },
\]
where we choose  $C$  large but independent of $\e, R_1$. Note that now $C$ is fixed.
\[
\D \F{\partial\psi_0}{\partial\bm n}-W_2\psi_0=\D\F{\e y_2 y_1^2 + (\e y_2 - 1)y_2^2}{r^3\ln^3 r}  + \F{1-\e y_2}{r\ln^3 r} + \frac{Cy_2}{r^{2}\ln^2 r}+ O\left(\frac1{r\ln^4 r}  \right).
\]

If $\e y_2 \ge \frac18, $ using $(y_1, y_2) \in \partial B_{\frac1\e}$ to get  $\frac{y_2}{r}\ge \frac1{20}$ for  $\e $ small.  From $0 \le y_2 \le \frac2{\e}, |y_1| \le r$,
\[
\left|  \D\F{\e y_2 y_1^2 + (\e y_2 - 1)y_2^2}{r^3\ln^3 r}  + \F{1-\e y_2}{r\ln^3 r} \right|\le \frac4{r\ln^3 r}.
\]
Hence
\[
\D \F{\partial\psi_0}{\partial\bm n}-W_2\psi_0 \ge \frac{C}{20r\ln^2 r} - \frac4{r\ln^3 r} +  O\left(\frac1{r\ln^4 r}  \right) \ge \frac{1}{r\ln^3 r}.
\]
If $ 0\le \e y_2 \le \frac18,$  then $\frac{y_2^2}{r^2} \le \frac{y_2^2}{y_1^2 + y_2^2} \le \frac{\e}2y_2 \le \frac1{16}$ for $(y_1, y_2) \in \partial B_{\frac1\e}$. Obviously $\D\F{\e y_2 y_1^2 }{r^3\ln^3 r} + \frac{Cy_2}{r^{2}\ln^2 r} \ge 0 $, so
\[
\D \F{\partial\psi_0}{\partial\bm n}-W_2\psi_0 \ge \D \F{1-\e y_2}{r\ln^3 r}\left(1-\frac{y_2^2}{r^2} \right) + O\left(\frac1{r\ln^4 r}  \right) \ge\frac78\frac{15}{16r\ln^3 r} +  O\left(\frac1{r\ln^4 r}  \right) \ge \frac{1}{2r\ln^3 r}
\]
where $R_1$ is large. Finally, letting $\psi=2\psi_0$, we get the lemma.
\end{proof}

\begin{lemma}\label{allortho}

There is an $\e_0 >0$, for any $\e \in (0, \e_0)$ and any solution $\phi$
of
\begin{equation}
\begin{cases}
-\Delta \phi - W_1\phi = f,  \qquad & \text{in } D_{\e} , \\
\D \F{\partial\phi}{\partial\bm n} - W_2 \phi = h,  & \text{on } \partial D_{\e}, \medskip\\
\D \int_{D_{\e}} \chi_\e Z_{i,\e}\phi=0, \qquad &i=0, 1.
\end{cases}
\end{equation}
It holds  that
\[
\|\phi\|_{L^\infty(D_{\e})} \le C(\|f\|_{**,  D_{\e}} + \|h\|_{*,  \partial D_{\e}}).
\]
where $C$ is independent of $\e$.

\end{lemma}
\begin{proof}
Choose $R_0=2R_1, R_1$ being the constant of Lemma \ref{barrier}. Then according to the properties of the barrier $\psi$ and Lemma \ref{sl}, we can finish the proof very similar to that of Lemma 4.2 in {\cite{DPM}}. Here we omit the details.
\end{proof}

Next we will establish an a priori estimate for solutions to the above equation under the  orthogonality condition with respect to $Z_{1,\e}$ only.

\begin{lemma}\label{pe}
For $\e$ sufficiently small, if $\phi$ solves
\begin{equation}\label{olp}
\begin{cases}
-\Delta \phi - W_1\phi = f,  \qquad & \text{in } D_{\e} , \\
\D \F{\partial\phi}{\partial\bm n} - W_2 \phi = h,  & \text{on } \partial D_{\e},\medskip \\
\D \int_{D_{\e}} \chi_\e Z_{1,\e}\phi=0,
\end{cases}
\end{equation}
then
\begin{equation}\label{apriori}
\|\phi\|_{L^\infty(D_{\e})} \le C(\|f\|_{**,  D_{\e}} + \|h\|_{*,  \partial D_{\e}}).
\end{equation}
where $C$ is independent of $\e$.

\end{lemma}

\begin{proof}
Let $\phi$ be the solution of (\ref{olp}).  Let $R>R_0+1$ be large and fixed. Set
\[
h_1(x)=\D\F{\ln \frac\delta\e-\ln |x|}{\ln \frac\delta\e-\ln R} \quad \text{and} \quad h_2(x)=\D\F{\frac1{|x|^{\sigma}}-(\frac{\e}{\delta})^{\sigma}}{\frac1{R^{\sigma}}-(\frac{\e}{\delta})^{\sigma}}\qquad \text{for} \quad R\le |x|\le \frac{\delta}{\e},
\]
and
\[
\hat{h}_1(y)=h_1(F_\e(y)), \qquad \hat{h}_2(y)=h_2(F_\e(y)),
\]
where $\sigma\in (0, 1)$  is  small and to be determined later. Obviously
\[
\Delta h_1\equiv 0, \qquad \Delta h_2(x)=\D\F{\sigma^2}{|x|^{\sigma+2}(\frac1{R^{\sigma}}-(\frac{\e}{\delta})^{\sigma})}.
\]
Let $\eta_0, \eta_1, \eta_2$ be radial non-negative smooth cut-off functions on $\mathbb{R}^2$ so that
\[
\eta_1(|x|)\equiv 1 \quad \text{in} \ B_R(0),  \quad \eta_1(|x|)\equiv 0 \quad \text{in} \ \mathbb{R}^2\setminus B_{R+1}(0), \qquad \eta_0(|x|)=\eta_1(|x|+2)
\]
 and
\[
\eta_2\equiv 1 \quad \text{in} \ B_{\frac{\delta}{4\e}}(0),  \quad \eta_2\equiv 0 \quad \text{in} \ \mathbb{R}^2\setminus B_{\frac{\delta}{3\e}}(0), \quad |\nabla \eta_2| \le \frac{C\e}{\delta}, \quad |\nabla^2 \eta_2| \le \frac{C\e^2}{\delta^2} \ \text{in} \  B_{\frac{\delta}{3\e}}(0)\setminus B_{\frac{\delta}{4\e}}(0),
\]
where $C$ is independent of $ \sigma, \delta, \e$.

Again we write
\[
\hat{\eta}_0(y)=\eta_0(F_\e (y)), \qquad \hat{\eta}_1(y)=\eta_1(F_\e(y)),\qquad \hat{\eta}_2(y)=\eta_2(F_\e(y)).
\]
and define
\[
\tilde{Z}_{0,\e}=\hat{\eta}_0Z_{0,\e}+(1-\hat{\eta}_0)\hat{\eta}_1 \hat{h}_1 Z_{0,\e} +  (1-\hat{\eta}_0)(1-\hat{\eta}_1)\hat{\eta}_2\hat{h}_1\hat{h}_2Z_{0,\e}.
\]
Given $\phi$, the solution to (\ref{olp}), let
\[
\tilde{\phi}=\phi + d_0\tilde{Z}_{0,\e}, \qquad \text{where} \ d_0=-\D\F{\int_{D_{\e}} Z_{0,\e}\chi_\e \phi}{\int_{D_{\e}} Z^2_{0,\e}\chi_\e }.
\]
Then estimate (\ref{apriori}) is a direct consequence of the coming claim.

{\bf  Claim.} Choose $R$ large and  fixed,  $\sigma$ small enough and then  fixed, there exists $\e_0$ such that for all $0 < \e < \e_0$,
\begin{equation}\label{9}
|d_0|\le C \left(\|f\|_{**,  D_{\e}}   +\|h\|_{*,  \partial D_{\e}}    \right).
\end{equation}

To prove the claim, observe, with the notation $L=-\Delta -W_1$, that
\begin{equation}\label{tildephi}
L(\tilde{\phi})=f+d_0L(\tilde{Z}_{0,\e}) \qquad \text{in} \ D_{\e}
\end{equation}
and
\[
\D \F{\partial\tilde{\phi}}{\partial\bm n} - W_2 \tilde{\phi}=h+d_0 \left(\D \F{\partial}{\partial\bm n} - W_2 \right)\tilde{Z}_{0, \e} \qquad \text{on} \ \partial D_{\e}.
\]
Then by Lemma \ref{allortho}, we have
\begin{equation}\label{phi}
\|\tilde{\phi}\|_{L^\infty(D_{\e})}\le C|d_0|\left(\left\| L(\tilde{Z}_{0,\e})  \right\|_{**, D_\e} + \left\| \left(\D \F{\partial}{\partial\bm n} - W_2 \right)\tilde{Z}_{0, \e}  \right\|_{*,\partial D_\e}\right)
+C\left(\|f\|_{**,  D_{\e}} + \|h\|_{*,  \partial D_{\e}}\right).
\end{equation}
Multiplying the equation (\ref{tildephi}) by $\tilde{Z}_{0,\e}$ and integrating by parts we find that
\begin{eqnarray*}
&&d_0 \left[\int_{D_{\e}} L(\tilde{Z}_{0,\e})\tilde{Z}_{0,\e} + \int_{\partial D_\e}\tilde{Z}_{0,\e} \left(\D \F{\partial}{\partial\bm n} - W_2 \right)\tilde{Z}_{0, \e}  \right]\\
=&&-\int_{D_\e} f\tilde{Z}_{0,\e} +  \int_{\partial D_\e}\tilde{\phi} \left(\D \F{\partial}{\partial\bm n} - W_2 \right)\tilde{Z}_{0, \e}-\int_{\partial D_\e} h\tilde{Z}_{0,\e}-\int_{D_\e} \tilde{\phi}L(\tilde{Z}_{0,\e})\\
\le && C\|f\|_{**,  D_\e}+C\|\tilde{\phi}\|_{L^\infty}\cdot\left(\left\|\F{\partial \tilde{Z}_{0, \e}}{\partial\bm n} - W_2 \tilde{Z}_{0, \e}\right\|_{*,  \partial D_\e}+\|L(\tilde{Z}_{0,\e})\|_{**,  D_\e}\right)+C\|h\|_{*, \partial D_\e}\\
\le && C\left( \|f\|_{**,  D_\e} + \|h\|_{*,  \partial D_\e}   \right)\left(\left\|\F{\partial \tilde{Z}_{0, \e}}{\partial\bm n} - W_2 \tilde{Z}_{0, \e}\right\|_{*,  \partial D_\e}+\|L(\tilde{Z}_{0,\e})\|_{**, D_\e}\right)\\
 && + C|d_0|\left(\left\|\F{\partial \tilde{Z}_{0, \e}}{\partial\bm n} - W_2 \tilde{Z}_{0, \e}\right\|^2_{*,  \partial D_\e}+\|L(\tilde{Z}_{0,\e})\|^2_{**,  D_\e}\right).
\end{eqnarray*}
due to (\ref{phi}) and the constant $C$ doesn't depend on $\e, \sigma$.\par

For $R$ large and fixed and for any $\sigma \in (0, 1)$ and  small $\e$, we may easily achieve the claim (\ref{9}) by building the ensuing estimates.
\begin{equation}\label{rd}
\int_{D_\e} L(\tilde{Z}_{0,\e})\tilde{Z}_{0,\e} + \int_{\partial D_\e}\tilde{Z}_{0,\e} \left(\D \F{\partial}{\partial\bm n} - W_2 \right)\tilde{Z}_{0, \e} \ge c_0\sigma -C\left( \frac1{\ln \frac1\e}+\e^\sigma\right) ,
\end{equation}
\begin{equation}\label{ld1}
\|L(\tilde{Z}_{0,\e})\|_{**,  D_\e} \le C \left(\sigma+\frac1{\ln \frac1\e}\right),
\end{equation}
\begin{equation}\label{ld2}
\left\|\F{\partial \tilde{Z}_{0, \e}}{\partial\bm n} - W_2 \tilde{Z}_{0, \e}\right\|_{*, \partial D_\e} \le C\left(\textcolor[rgb]{1.00,0.00,0.00}{\sigma} + \frac1{\ln \frac1\e}\right),
\end{equation}
where the positive constants $c_0, C$ are independent of $\sigma, \e$.

\medskip
{\bf  Proof of (\ref{rd}).} Denote $x=F_\e(y)$ hereinafter and recall that this map preserves area. It is divided that
\[\int_{ D_\e} L(\tilde{Z}_{0,\e})\tilde{Z}_{0,\e} = I_1 + I_2 + I_3+I_4 + I_5+I_6:=\sum_{i=1}^6\int_{B_i}  L(\tilde{Z}_{0,\e})\tilde{Z}_{0,\e},\]
where
\[
B_1 = (F_\e)^{-1}(\{r<R-2\}\cap \mathbb{R}^2_+), \qquad \quad  B_2=(F_\e)^{-1}(\{R-2<r<R-1\}\cap \mathbb{R}^2_+),
\]
\[
 B_3=(F_\e)^{-1}(\{R-1<r<R\}\cap \mathbb{R}^2_+), \qquad B_4=(F_\e)^{-1}(\{R<r<R+1\}\cap \mathbb{R}^2_+),
\]
and
\[
B_5=(F_\e)^{-1}\left(\left\{R+1<r<\frac{\delta}{4\e}\right\}\cap \mathbb{R}^2_+\right),  \quad B_6=(F_\e)^{-1}\left(\left\{\frac{\delta}{4\e}<r<\frac{\delta}{3\e}\right\}\cap \mathbb{R}^2_+\right) \ \ \text{with} \  r=|x|.
\]
Then we estimate term by term. In $x$, the operator $L$ has the form that
\begin{align*}
\tilde{L}&= -\Delta  + O(\e|x|)\nabla^2 + O(\e)\nabla - W_1\left((F_\e)^{-1}x\right)\\
&= -\Delta  + O(\e|x|)\nabla^2 + O(\e)\nabla - W_1\left(x+O(\e |x|)\right)\\
&= -\Delta-\frac{8\lambda}{\lambda^2 + x_1^2 + (x_2 +\frac{b}{\sqrt{a}}\lambda )^2} + O(\e|x|)\nabla^2 + O(\e)\nabla+O\left(\frac{\e^\alpha + \e|x|}{1+|x|^4}\right)
\end{align*}
with $a=K(\xi), b=\kappa(\xi)$. \par

\emph{Estimate of $I_1$.} In $B_1$, $\tilde{Z}_{0, \e} =Z_{0, \e}$. So we get that
\begin{eqnarray*}
I_1 &=& \int_{\{0 \le r \le R-2\}\cap \mathbb{R}^2_+} \tilde{L}(z_0)z_0\\
& =&\int_{\{0\le r \le R-2\}\cap \mathbb{R}^2_+} \left(-\Delta z_0  + O(\e|x|)\nabla^2 z_0 + O(\e)\nabla z_0 - W_1(\xi'+x+O(\e |x|))z_0\right)z_0\\
&=& O\left( \e^\alpha  + \e \ln R \right)=O(\e^{\alpha}).
\end{eqnarray*}

\emph{Estimate of $I_2$.} In $B_2$, $\tilde{Z}_{0, \e} =Z_{0, \e} + (1-\hat{\eta}_0)(\hat{h}_1-1)Z_{0,\e}$.  Let $\tilde{z}_0=z_0 + (1-\eta_0)(h_1 -1)z_0$. Then it holds that
\begin{align*}
I_2 =&\ \int_{\{R-2 \le r \le R-1\}\cap \mathbb{R}^2_+} \tilde{L}(\tilde{z}_0)\tilde{z}_0\\
=&\ \int_{\{R-2\le r \le R-1\}\cap \mathbb{R}^2_+} \Big\{-z_0\Delta [(1-\eta_0)(h_1-1)] -2\nabla z_0\cdot \nabla  [(1-\eta_0)(h_1-1)] + O(\e|x|)\nabla^2\tilde{z}_0  \\
&\qquad\qquad\qquad\qquad\qquad + O(\e)\nabla \tilde{z}_0 \Big\}\tilde{z}_0+O(\e^{\alpha})\\
=&\ O\left( \e^\alpha  + \frac1{\ln \frac\delta\e} \right).
\end{align*}\par

\emph{Estimate of $I_3$.}  In $B_3$,  $\tilde{Z}_{0, \e} =\hat{h}_1Z_{0,\e}$. We arrive at
\begin{eqnarray*}
I_3 &=& \int_{\{R-1 \le r \le R\}\cap \mathbb{R}^2_+} \tilde{L}(h_1z_0)h_1z_0\\
&=&\int_{\{R-1\le r \le R\}\cap \mathbb{R}^2_+} \left(-2\nabla z_0\cdot \nabla h_1  + O(\e|x|)\nabla^2(h_1z_0) + O(\e)\nabla (h_1z_0) \right)h_1z_0+O(\e^{\alpha})\\
&=& O\left( \e^\alpha  + \frac1{\ln \frac\delta\e} \right).
\end{eqnarray*}\par

\emph{Estimate of $I_4$.} Now $ \tilde{Z}_{0, \e} =\hat{\eta}_1\hat{h}_1Z_{0, \e} + (1-\hat{\eta}_1)\hat{h}_2\hat{h}_1Z_{0,\e}$
and set again $\tilde{z}_0(x) = h_1z_0(x) + (1-\eta_1)(h_2-1)h_1z_0(x).$ Then
\begin{eqnarray*}
I_4 &=& \int_{\{R \le r \le R+1\}\cap \mathbb{R}^2_+} \tilde{L}(\tilde{z}_0)\tilde{z}_0\\
&=& \int_{\{R < r < R+1\}\cap \mathbb{R}^2_+} \left\{ -2\nabla h_1\cdot \nabla z_0-\Delta[(1-\eta_1)(h_2-1)h_1]z_0 -2\nabla  [(1-\eta_1)(h_2-1)h_1]\cdot \nabla z_0 \right\}\tilde{z}_0+O(\e^{\alpha})\\
&=& \int_{\{R < r < R+1\}\cap \mathbb{R}^2_+}  -\Delta[(1-\eta_1)(h_2-1)h_1]z_0\tilde{z}_0 +O\left(\frac\sigma{R^2}+\frac1{\ln \frac\delta\e}+\e^{\alpha}\right).
\end{eqnarray*}
The functions  $\eta_1$, $h_2$ are both radial symmetric and decreasing with $ h_2'(r) = \D\F{-\sigma}{|r|^{\sigma+1}\left(\frac1{R^{\sigma}}-\frac{\e^{\sigma}}{\delta^{\sigma}}    \right)}$, $\eta_1'(R)=\eta_1(R+1)=\eta_1'(R+1)=0$, $\eta_1(R)=1$, $\left.\frac{\partial \tilde{h}}{\partial x_2}\right|_{x_2=0}=\left.\frac{\partial \eta_1}{\partial x_2}\right|_{x_2=0}=0$. This leads to
\begin{eqnarray*}
&&\int_{\{R < r < R+1\}\cap \mathbb{R}^2_+}  -\Delta[(1-\eta_1)(h_2-1)h_1]z_0\left[h_1z_0(x) + (1-\eta_1)(h_2-1)h_1z_0(x)\right]\\
=&& -\int_{\partial\left(\{R < r < R+1\}\cap \mathbb{R}^2_+\right)} \frac{\partial [(1-\eta_1)(h_2-1)h_1]}{\partial\nu} \left[h_1z^2_0(x) + (1-\eta_1)(h_2-1)h_1z^2_0(x)\right]\\
&& + \int_{\{R < r < R+1\}\cap \mathbb{R}^2_+} \nabla [(1-\eta_1)(h_2-1)h_1] \cdot \nabla \left[h_1z^2_0(x) + (1-\eta_1)(h_2-1)h_1z^2_0(x)\right]\\
=&& -\int_{\{r=R+1\}\cap \mathbb{R}^2_+} h_2'(r)h_1^2h_2z^2_0(x)+  O\left(\frac{\sigma^2}{R}+\frac\sigma{R^2}+\frac{1}{\ln \frac\delta\e}\right)
\end{eqnarray*}
where $\nu$ denotes the unit outer normal vector.\par

\emph{Estimate of $I_5$.} In $B_5$, $\tilde{Z}_{0,\e} = \hat{h}_1\hat{h}_2Z_{0,\e}$ and set $\tilde{z}_0=h_1h_2z_0$. Then
\begin{align}
I_5=&\ \int_{\{R+1 <r <\frac\delta{4\e}\}\cap \mathbb{R}^2_+} \tilde{L}(\tilde{z}_0)\tilde{z}_0 \nonumber\\
=&\ \int_{\{R+1 <r <\frac\delta{4\e}\}\cap \mathbb{R}^2_+}\Big[-\Delta(h_1h_2) z_0 -2\nabla (h_1h_2)\cdot \nabla z_0 +  O(\e|x|)\nabla^2(h_1h_2z_0) + O(\e)\nabla (h_1h_2z_0) \Big]h_1h_2z_0 + O(\e^\alpha)\nonumber\\
=&\  -\int_{\partial(\{R+1 <r <\frac\delta{4\e}\}\cap \mathbb{R}^2_+)} \frac{\partial (h_1h_2)}{\partial \nu} h_1h_2z_0^2+O\left(\e^\alpha+\e^\sigma R^\sigma\sigma^2 + \frac1{\ln\frac\delta\e}\right)\nonumber\\
&\qquad\qquad +\int_{\{R+1 <r <\frac\delta{4\e}\}\cap \mathbb{R}^2_+}\left[ \nabla(h_1h_2) \cdot \nabla (h_1h_2z_0^2) - 2\nabla (h_1h_2)\cdot \nabla z_0 h_1h_2z_0 \right]\nonumber\\
=&\ \int_{\{r=R+1\}\cap \mathbb{R}^2_+} (h_1h_2)'(r)h_1h_2z_0^2 + \int_{\{R+1 <r <\frac\delta{4\e}\}\cap \mathbb{R}^2_+} |\nabla (h_1h_2)|^2z_0^2 +O\left(\e^\alpha+\e^\sigma R^\sigma + \frac1{\ln\frac\delta\e}\right)   \nonumber\\
=&\  \int_{\{r=R+1\}\cap \mathbb{R}^2_+} h_2'(r)h^2_1h_2z_0^2 + \int_{\{R+1 <r <\frac\delta{4\e}\}\cap \mathbb{R}^2_+}|\nabla h_2|^2h_1^2z_0^2 +  O\left(\e^\alpha+\e^\sigma R^\sigma + \frac1{\ln\frac\delta\e}\right)  \nonumber\\
=&\  \int_{\{r=R+1\}\cap \mathbb{R}^2_+} h_2'(r)h^2_1h_2z_0^2 + \int_{\{R+1 <r <\frac\delta{4\e}\}\cap \mathbb{R}^2_+} \frac{\sigma^2R^{2\sigma}}{|x|^{2+2\sigma}} z_0^2 +O\left(\e^\alpha+\e^\sigma R^\sigma + \frac1{\ln\frac\delta\e}\right)
\end{align}
and
\[
\int_{\{R+1 <r <\frac\delta{4\e}\}\cap \mathbb{R}^2_+} \frac{\sigma^2R^{2\sigma}}{|x|^{2+2\sigma}} z_0^2 \ge \gamma_0\sigma
\]
where the positive constant $\gamma_0$ is independent of $R, \e, \sigma$.

\emph{Estimate of $I_6$.} In $B_6$, $\tilde{Z}_{0,\e}=\hat{\eta}_2\hat{h}_1\hat{h}_2Z_{0,\e}$ and set $\tilde{z}_0=\eta_2h_1h_2z_0$. Now $h_1(r)=O\left(\frac1{\ln \frac\delta{\e}}\right), h_2(r)=R^\sigma\e^\sigma$ for $\frac\delta{4\e}<r<\frac\delta{3\e}$. Thus
\begin{eqnarray*}
I_6 &=& \int_{\{\frac\delta{4\e} \le r \le \frac\delta{3\e}\}\cap \mathbb{R}^2_+} \tilde{L}(\eta_2h_1h_2z_0)\eta_2h_1h_2z_0\\
&=&\int_{\{\frac\delta{4\e} \le r \le \frac\delta{3\e}\}\cap \mathbb{R}^2_+} \left(-\Delta(\eta_2h_1h_2)z_0-2\nabla z_0\cdot \nabla (\eta_2h_1h_2)  + O(\e|x|)\nabla^2\tilde{z}_0 + O(\e)\nabla \tilde{z}_0 \right)\eta_2h_2h_1z_0+O(\e^{\alpha})\\
&=& O\left( \e^\alpha  + \frac1{\ln \frac\delta\e} \right).
\end{eqnarray*}\par

From the above six estimates, we conclude that
\BE\label{8}
\int_{ D_\e} L(\tilde{Z}_{0,\e})\tilde{Z}_{0,\e} \geq \gamma_0\sigma + O\left( \e^\alpha  + \frac1{\ln \frac\delta\e} +\frac{\sigma^2}{R}+\frac\sigma{R^2}\right) .
\EE \par

Next we estimate  $\int_{\partial D_\e}\tilde{Z}_{0,\e} \left(\D \F{\partial}{\partial\bm n} - W_2 \right)\tilde{Z}_{0, \e}$. It is checked, writing $x=(x_1,x_2)$, that
\[
\int_{\partial D_\e}\tilde{Z}_{0,\e} \left(\D \F{\partial}{\partial\bm n} - W_2 \right)\tilde{Z}_{0, \e} =\int_{\partial \mathbb{R}^2_+ } \tilde{z}_0 \Big[B(\tilde{z}_0) - W_2(\xi'+x_1+O(\e |x_1|) \tilde{z}_0\Big]j(x)
\]
where $\tilde{z}_0=\tilde{Z}_{0,\e}(F_\e^{-1}(x))$ and $j(x)$ is a positive function arising from the change of variables bounded uniformly in  $\e$. $B$ is  a differential operator of order one on $\partial \mathbb{R}^2_+$. Rotating $\partial D_\e$ so that $\nabla F_\e(\xi')=I$ we get the following expansion
\[
B=-\frac{\partial}{\partial x_2} + O(\e |x_1|)\nabla.
\]
Recall that
\begin{eqnarray*}
W_2(\xi'+x_1+O(\e |x_1|)&=&\frac{\kappa(\e y)}{\sqrt{K(\xi)}}\D\F{2\lambda}{\lambda^2+|y-\xi'-\mathfrak D(\xi)\lambda\bm n(\xi)|^2}\left[1+O(\e|y-\xi'|)+O(\e^\alpha)\right]\\
&=&\frac{2b\lambda}{\sqrt{a} \left(\lambda^2+x_1^2 + \frac{b^2\lambda^2}{a}\right)} \left[1+ O\left( \e|x_1|+\e^\alpha\right) \right].
\end{eqnarray*}
Thus
\[
\int_{\partial \mathbb{R}^2_+ \cap\{ |x_1|\le R-2\}} z_0 \Big[B(z_0) - W_2(\xi'+x_1+O(\e |x_1|) z_0\Big]j(x)=O(\e^\alpha +\e\ln R).
\]
Note that $\tilde{z}_0= z_0 +(1-\eta_0)(h_1-1)z_0$  for $R-2 \le |x|\le R-1$, using the fact that $\eta_0, h_1$ has zero normal derivative on $\mathbb{R}^2_+ $, we have
\[
B(\tilde{z}_0)=-\frac{\partial z_0}{\partial x_2} - (1-\eta_0)(h_1-1)\frac{\partial z_0}{\partial x_2} + O(\e |x_1|)\nabla \tilde{z}_0 \quad  \text{with} \ \quad x_2=0.
\]
It is easy to see that
\[
\int_{\partial \mathbb{R}^2_+ \cap \{R-2\le |x_1|\le R-1\}}\tilde{z}_0 \Big[B(\tilde{z}_0) - W_2(\xi'+x_1+O(\e |x_1|) ) \tilde{z}_0\Big]j(x)=O\left(\e^\alpha +\frac1{\ln \frac\delta{\e}}\right).
\]
Similarly,
\[
\int_{\partial \mathbb{R}^2_+ \cap \{R-1\le |x_1|\le R\}}h_1z_0 \Big[B(h_1z_0) - W_2(\xi'+x_1+O(\e |x_1|)) h_1z_0\Big]j(x)=O\left(\e^\alpha\right)
\]
and
\[
\int_{\partial \mathbb{R}^2_+ \cap \{R\le |x_1|\le R+1\}} \tilde{z}_0 \Big[B(\tilde{z}_0) - W_2(\xi'+x_1+O(\e |x_1|)) \tilde{z}_0\Big]j(x)=O\left(\e^\alpha\right)
\]
with $\tilde{z}_0 =\eta_1h_1z_0+(1-\eta_1)h_1h_2z_0$. Also we obtain that
\begin{multline*}
\int_{\partial \mathbb{R}^2_+ \cap \{R+1\le |x_1|\le \frac\delta{4\e}\}} h_1h_2z_0 \Big[B(h_1h_2z_0) - W_2(\xi'+x_1+O(\e |x_1|)) h_1h_2z_0\Big]j(x)\\
=O\left(\e^\alpha+\e\ln\frac1\e + \frac1{\ln \frac1\e} + \e^\sigma R^\sigma\right),
\end{multline*}
and
\begin{multline*}
\int_{\partial \mathbb{R}^2_+ \cap \{\frac\delta{4\e}\le |x_1|\le \frac\delta{3\e}\}} \eta_2h_1h_2z_0 \Big[B(\eta_2h_1h_2z_0) - W_2(\xi'+x_1+O(\e |x_1|)) \eta_2h_1h_2z_0\Big]j(x) \\= O\left(\e^\alpha+R^{2\sigma}\e^{2\sigma}\ln^{-2}\frac1\e\right).
\end{multline*}

Finally combing the above estimate and (\ref{8}), we have that
\BEN
\int_{D_\e} L(\tilde{Z}_{0,\e})\tilde{Z}_{0,\e} + \int_{\partial D_\e}\tilde{Z}_{0,\e} \left(\D \F{\partial}{\partial\bm n} - W_2 \right)\tilde{Z}_{0, \e} \geq
 \gamma_0\sigma+O\left(\frac{\sigma^2}{R} + \frac{\sigma}{R^2} + \frac1{\ln \frac1\e} + \e^\sigma R^\sigma\right).
\EEN
Now  take $R$ large enough, we know (\ref{rd}) holds.
\medskip\par

{\bf Proof of (\ref{ld1}).}
For $y \in B_1$, $\tilde Z_{0,\e} = Z_{0,\e}$. Obviously we have
\begin{multline*}
L(\tilde{Z}_{0,\e})=\tilde{L}(z_0)=-\Delta z_0 + O(\e|x|)\nabla^2z_0 + O(\e)\nabla z_0 - W_1\left((\xi' + x + O(\e |x|)\right)z_0\\
= O(\e|x|)\nabla^2z_0 + O(\e)\nabla z_0 + O(\e |x|+ \e^\alpha)z_0=O(\e |x|+ \e^\alpha),
\end{multline*}
and hence
\[
\|L(\tilde{Z}_{0,\e})\|_{**,  B_1} \le C  \e^{\alpha}.
\]\par

For $y \in B_2$, then $\tilde{z}_0=z_0 + (1-\eta_0)(h_1-1)z_0$. So
\begin{align*}
L(\tilde{Z}_{0,\e})&=\tilde{L}(\tilde{z}_0)=-\Delta \tilde{z}_0 +  O(\e|x|)\nabla^2 \tilde{z}_0 + O(\e)\nabla \tilde{z}_0 - W_1\left((\xi' + x + O(\e |x|)\right)\tilde{z}_0\\
&= -z_0\Delta[(1-\eta_0)(h_1-1)]-2\nabla z_0\cdot\nabla [(1-\eta_0)(h_1-1)] + O(\e^\alpha)\\
&= O\left( \e^\alpha + \frac1{\ln \frac1\e}\right),
\end{align*}
which leads to
\[
\|L(\tilde{Z}_{0,\e})\|_{**,  B_2} \le  \frac C{\ln \frac1\e}.
\]\par

For $y \in B_3$, it holds that $\tilde{z}_0=h_1z_0$. Then
\begin{eqnarray*}
L(\tilde{Z}_{0,\e})&=&\tilde{L}(\tilde{z}_0)=-h_1\Delta z_0-2\nabla h_1 \cdot \nabla z_0 +  O(\e|x|)\nabla^2 \tilde{z}_0 + O(\e)\nabla \tilde{z}_0 - W_1\left((\xi' + x + O(\e |x|)\right)\tilde{z}_0\\
&=& O\left( \e^\alpha + \frac1{\ln \frac1\e}\right),
\end{eqnarray*}
from which we can get that
\[
\|L(\tilde{Z}_{0,\e})\|_{**,  B_3} \le  \frac C{\ln \frac1\e}.
\]\par

For $y \in B_4$, obviously $\tilde{z}_0=h_1z_0 + (1-\eta_1)(h_2-1)h_1z_0$. Note that $h_2-1 = \F{\F{1}{|x|^\sigma}-\F{1}{R^\sigma}}{\F{1}{R^\sigma}-\e^\sigma} \leq C\sigma$ in this region. We have
\begin{align*}
L(\tilde{Z}_{0,\e})&=\tilde{L}(\tilde{z}_0)=-2\nabla h_1\cdot \nabla z_0-z_0\Delta[(1-\eta_1)(h_2-1)h_1] - 2\nabla z_0\cdot \nabla [(1-\eta_1)(h_2-1)h_1]+ O(\e^\alpha)\\
&= O\left( \e^\alpha + \frac1{\ln \frac1\e} + \sigma\right),
\end{align*}
from which we can get that
\[
\|L(\tilde{Z}_{0,\e})\|_{**,  B_4} \le  \frac C {\ln \frac1\e} + C\sigma.
\]\par

For $y \in B_5$, we see that $\tilde{z}_0=h_1h_2z_0$. Thus it is checked that
\begin{align*}
&\ O(\e |x|) \nabla^2 (h_1 h_2 z_0) \\
=&\ O(\e |x|)\left[ \nabla^2h_1 h_2 z_0 + h_1 \nabla^2 h_2 z_0 +  h_1  h_2 \nabla^2 z_0 +2 \nabla h_1 \nabla h_2 z_0 + 2 \nabla h_1 h_2 \nabla z_0 +2h_1 \nabla h_2 \nabla z_0  \right]\\
=&\ O\left( \F{\e}{|x|^{1+\sigma}} \right),
\end{align*}
\[
O(\e)\nabla (h_1h_2z_0) =  O\left(\frac \e {\ln \frac1\e |x|}\left[\frac1{|x|^\sigma} - (\frac\e\delta)\sigma\right] + \frac{\ln\frac\delta\e-\ln |x|}{\ln\frac\delta\e}\frac{\e\sigma}{|x|^{1+\sigma}} + \frac \e {|x|^{2+\sigma}}    \right) = O\left( \F{\e}{|x|^{1+\sigma}} \right),
\]
and
\[
-\Delta (h_1h_2z_0) - W_1\left((\xi' + x + O(\e |x|)\right)h_1h_2z_0= O\left( \frac{\sigma}{\ln \frac1\e |x|^{2+\sigma}}+\frac{\sigma^2}{|x|^{2+\sigma}}+\frac{1}{\ln \frac1\e |x|^{3}}+\frac\sigma{|x|^{3+\sigma}}+\frac{\e|x| + \e^\alpha}{1+|x|^4}\right).
\]
It is easy to get that
\[
\|L(\tilde{Z}_{0,\e})\|_{**,  B_5} \le C\sigma + C\e^{\sigma}\ln^3\frac1\e+C\e^{\alpha}+ \frac C{\ln \frac1\e}.
\]\par

For $y \in B_6$, $\tilde{z}_0 = \eta_2h_1h_2z_0$, $h_1=O(\frac1{\ln \frac1\e})$, $h_2=O(\e^\sigma)$.  Hence
\begin{align*}
L(\tilde{Z}_{0,\e})&=\tilde{L}(\tilde{z}_0)=-\Delta \tilde{z}_0 +  O(\e|x|)\nabla^2 \tilde{z}_0 + O(\e)\nabla \tilde{z}_0 - W_1\left((\xi' + x + O(\e |x|)\right)\tilde{z}_0\\
&= O\left( \frac{\e^{2+\sigma}}{\ln\frac1\e} +  \frac{\e|x| + \e^\alpha}{1+|x|^4}\right) .
\end{align*}
and
\[
\|L(\tilde{Z}_{0,\e})\|_{**,  B_6} \le C\e^\sigma\ln^2\frac1\e + C\e^4.
\]\par

From the above estimates, (\ref{ld1}) is finally concluded.

\medskip
{\bf Proof of (\ref{ld2}).} It is known, in the proof of (\ref{rd}), that $\F{\partial \tilde{Z}_{0, \e}}{\partial\bm n} - W_2 \tilde{Z}_{0, \e}$ is transformed to
$B(\tilde{z}_0) - W_2(\xi'+x_1+O(\e |x_1|) \tilde{z}_0$. We still estimate it in $B_i$ respectively.
For $y \in B_1$, $\tilde z_0 = z_0$.
\[
B(z_0)- W_2(\xi'+x_1+O(\e |x_1|) z_0 = O(\e^\alpha).
\]
For $y \in B_2$, $\tilde{z}_0=z_0 + (1-\eta_0)(h_1-1)z_0$.
\[
B(\tilde{z}_0)- W_2(\xi'+x_1+O(\e |x_1|) \tilde{z}_0=O\left( \e^\alpha + \frac1{\ln\frac1\e} \right).
\]
For $y \in B_3$, $\tilde{z}_0=h_1z_0$,
\[
B(h_1z_0)- W_2(\xi'+x_1+O(\e |x_1|) h_1 z_0 = O(\e^\alpha).
\]
For $y \in B_4$, $\tilde{z}_0=h_1z_0 + (1-\eta_1)(h_2-1)h_1z_0$
\[
B(\tilde{z}_0)- W_2(\xi'+x_1+O(\e |x_1|) \tilde{z}_0=O(\e^\alpha).
\]
For $y \in B_5$, $\tilde{z}_0=h_1h_2z_0,$
\[
B(h_1h_2z_0)- W_2(\xi'+x_1+O(\e |x_1|) h_1h_2 z_0 = O\left(\frac\e{|x_1|} + \frac\e{\ln\frac1\e |x_1|^\sigma} + \frac{\e|x_1|+\e^\alpha}{1+x_1^2}\right).
\]
And  for $y \in B_6, \tilde{z}_0 = \eta_2h_1h_2z_0, h_1=O(\frac1{\ln \frac1\e}), \ h_2=O(\e^\sigma)$,
\[
B(\eta_2h_1h_2z_0)- W_2(\xi'+x_1+O(\e |x_1|) \eta_2h_1h_2 z_0 = O\left(\frac{\e^{1+\sigma}}{\ln\frac1\e}+ \frac{\e^\sigma}{|x_1|\ln\frac1\e}+ \frac{\e|x_1|+\e^\alpha}{1+x_1^2}\right).
\]
Combining all these estimates and the definition of $\|\cdot\|_{*,\partial D_\e}$, we get the desired (\ref{ld2}).
\end{proof}

\begin{proposition}\label{lpe}
There exist an  $\e_0>0$ and a positive number $C$ such that for any $\xi' \in \partial D_\e$ and $f \in L^\infty (D_\e), h\in L^\infty(\partial D_\e)$, there is a unique solution $\phi \in L^\infty(D_\e), c_1 \in \mathbb{R}$ to
\begin{equation}\label{lp}
\begin{cases}
-\Delta \phi - 2K(\e y)e^{2V}\phi = f  + c_1\chi_\e Z_{1,\e},  \qquad & \text{in } D_\e , \\
\D \F{\partial\phi}{\partial\bm n} - \kappa(\e y) e^V \phi = h, \quad & \text{on } \partial D_\e, \medskip \\
\D\int_{D_\e} \chi_\e  Z_{1,\e}\phi=0.
\end{cases}
\end{equation}
Moreover,
\[
\|\phi\|_{L^\infty(D_\e)} \le C\|f\|_{**, D_\e} +C\|h\|_{*,\partial D_\e}.
\]

\end{proposition}

\begin{proof}
By Lemma \ref{pe}, any solution to (\ref{lp}) satisfies
\[
\|\phi\|_{L^\infty(D_\e)} \le C(\|f\|_{**,  D_\e} +|c_1|+ \|h\|_{*,  \partial D_\e});
\]
Therefore, for the estimate, it is enough to prove
\[
|c_1| \le C(\|f\|_{**,  D_\e} + \|h\|_{*,  \partial D_\e}).
\]

In the following, let $\hat{\eta}_2$ be defined as before.  Multiplying the first equation of (\ref{lp}) by $\hat{\eta}_2Z_{1,\e}$ and integrating by parts we can get that
\begin{multline*}
c_1\int_{D_\e} \chi_\e Z^2_{1,\e}=-\int_{D_\e} f\hat{\eta}_2Z_{1,\e} + \int_{D_\e}\left[-\Delta(\hat{\eta}_2Z_{1, \e})-W_1\hat{\eta}_2Z_{1,\e}\right]\phi\\
-\int_{\partial D_\e} h\hat{\eta}_2Z_{1,\e} + \int_{\partial D_\e} \phi\frac{\partial \hat{\eta}_2}{\partial \bm n}Z_{1,\e} + \int_{\partial D_\e} \phi\hat{\eta}_2\left(\frac{\partial Z_{1,\e}}{\partial \bm n}-W_2Z_{1,\e}\right).
\end{multline*}
From $|\nabla \hat{\eta}_2| \le C\e,  |\nabla^2 \hat{\eta}_2|\le C\e^2,  Z_{1, \e} =O(\frac1{1 + |y-\xi'|})$, direct computation leads to
\[
\left| -\int_{D_\e} f\hat{\eta}_2Z_{1,\e} -\int_{\partial D_\e} h\hat{\eta}_2Z_{1,\e}  \right| \le C(\|f\|_{**,  D_\e} + \|h\|_{*,  \partial D_\e}),
\]

\[
\int_{D_\e}\left[-\Delta(\hat{\eta}_2Z_{1, \e})-W_1\hat{\eta}_2Z_{1,\e}\right]\phi =O\left( \e^\alpha + \e\ln\frac1\e\right)\|\phi\|_{L^\infty(D_\e)},
\]
and
\[
\int_{\partial D_\e} \phi\frac{\partial \hat{\eta}_2}{\partial \bm n}Z_{1,\e} + \int_{\partial D_\e} \phi\hat{\eta}_2\left(\frac{\partial Z_{1,\e}}{\partial \bm n}-W_2Z_{1,\e}\right) =O\left( \e^\alpha + \e\ln\frac1\e\right)\|\phi\|_{L^\infty(D_\e)}.
\]
Then it is easy to obtain
\[
|c_1|\le C(\|f\|_{**,  D_\e} + \|h\|_{*,  \partial D_\e}) + C\left( \e^\alpha + \e\ln\frac1\e\right)\|\phi\|_{L^\infty(D_\e)},
\]
from which we can get
\begin{equation}
\|\phi\|_{L^\infty(D_\e)} \le C(\|f\|_{**,  D_\e} + \|h\|_{*,  \partial D_\e}).
\end{equation}

Now consider the Hilbert space
\[
H=\left\{ \phi \in H^1(D_\e): \int_{D_\e} \chi_\e Z_{1,\e}\phi=0\right\}
\]
with the norm $\|\phi\|^2_{H}=\int_{D_\e}  |\nabla \phi|^2 + |\phi|^2$. Equation (\ref{lp}) is equivalent to find $\phi \in H$ such that
\[
\int_{D_\e} \left(\nabla\phi\cdot\nabla\psi + \phi\psi \right) - \int_{D_\e} \left(\phi\psi+W_1\phi\psi\right)-\int_{\partial D_\e} W_2\phi\psi=\int_{D_\e} f\psi + \int_{\partial D_\e} h\psi, \quad \forall \psi \in H.
\]
By Fredholm's alternative it is enough to show the uniqueness of solutions to the  problem (\ref{lp}), which is guaranteed by Lemma \ref{pe}.
\end{proof}

\BR
The result of Proposition \ref{lpe} implies that the unique solution $\phi=T(f,h)$ of (\ref{lp}) defines a continuous linear map. For later purposes, the differentiability of the operator $T$ with respect to $\xi'$ is necessary. The proof is almost the same as that in \cite{DPM}. Here we omit the details.
\ER

\section{The nonlinear problem}\label{s4}

In this section, we will solve the following nonlinear problem that
\begin{equation} \label{nonlp}
\begin{cases}
-\Delta \phi - 2K(\e y)e^{2V}\phi =R_1(y) + N_1(\phi)+ c_1\chi_\e Z_{1,\e}, \qquad \qquad  & \text{in } D_{\e} , \\
\D \F{\partial\phi}{\partial\bm n} - \kappa(\e y) e^V \phi =R_2(y) + N_2(\phi) , & \text{on }\partial D_{\e}, \medskip \\
\D \int_{D_\e} \chi_\e Z_{1,\e} \phi=0,
\end{cases}
\end{equation}
where
\[
N_1(\phi)= K(\e y) e^{2V}(e^{2\phi}-1-2\phi), \qquad  N_2(\phi)=\kappa(\e y) e^{V}(e^{\phi}-1-\phi).
\]
Recall in Lemma \ref{lem1} that
\begin{align*}
  &\ R_1(y)  = \F{4\lambda^2}{(\lambda^2+|y-\xi'-\mathfrak D(\xi)\lambda\bm n(\xi)|^2)^2} \left[O(\e |y-\xi'|)+O(\e^\alpha)\right], \\
  \\
  &\ R_2(y)=\F{\mathfrak D(\xi) 2\lambda}{(\lambda^2+|y-\xi'-\mathfrak D(\xi)\lambda\bm n(\xi)|^2)} \left[O(\e |y-\xi'|)+O(\e^\alpha)\right] + \e \F{d\lambda^2}{\lambda^2 + |y-\xi'|^2},
\end{align*}
in $|y-\xi'|\leq \F{\delta}{\e}$,
from which we may obtain that
\begin{equation}\label{error}
\|R_1(y)\|_{**, D_\e} \le C\e^\alpha,  \qquad \qquad \|R_2(y)\|_{*, \partial D_\e} \le C\e^\alpha, \qquad \text{for any } \ \alpha \in (0, 1).
\end{equation}

\begin{lemma}\label{nonlinear}
There exists $\e_0>0$ and $C>0$, such that for $0< \e <\e_0$ and $\xi' \in \partial D_\e$, the problem (\ref{nonlp}) admits a unique solution $(\phi, c_1)$ such that
\begin{equation}
\|\phi\|_{L^\infty(D_\e)} + |c_1| \le C\e^\alpha, \qquad \qquad \text{for any} \ \ \alpha \in (0, 1).
\end{equation}
Furthermore, the function $\xi'\rightarrow \phi(\xi') \in C(\bar{D}_{\e} )$ is $C^1$ and
\begin{equation}
\|D_{\xi'} \phi\|_{L^\infty(D_\e)}  \le C\e^\alpha.
\end{equation}
\end{lemma}

\begin{proof}
In terms of the operator $T$ defined in the previous section, problem (\ref{nonlp}) becomes
\begin{equation}
\phi=T(R_1+N_1(\phi), R_2+N_2(\phi))\equiv A(\phi).
\end{equation}
For a given number $\gamma>0$, let us consider the region
\[
\mathcal{F}_\gamma =\left\{\phi \in H: \|\phi\|_{L^\infty(D_\e)} \le \gamma \e^{\alpha}\right\}.
\]
From Proposition \ref{lpe}, we get
\begin{align*}
\|A(\phi)\|_{L^\infty(D_\e)} &\le  C\left(\|R_1+N_1(\phi)\|_{**, D_\e} + \|R_2+N_2(\phi)\|_{*, \partial D_\e}   \right)\\
&\le  C\left(\|R_1\|_{**, D_\e} + \|R_2\|_{*, \partial D_\e} +\|N_1(\phi)\|_{**, D_\e} + \|N_2(\phi)\|_{*, \partial D_\e}\right)\\
&\le C\e^\alpha + C\left\|\frac{1}{1+|y-\xi'|^4}\phi^2\right\|_{**, D_\e}+ C\left\|\frac{1}{1+|y-\xi'|^2}\phi^2\right\|_{*, \partial D_\e}\\
&\le C\e^\alpha + \gamma^2 \e^{2\alpha}
\end{align*}
which leads to $A(\phi) \in  \mathcal{F}_\gamma$  for any $\phi \in \mathcal{F}_\gamma$ where $\gamma$ is large but fixed. Also, for any $\phi_1, \phi_2 \in \mathcal{F}_\gamma$,
\[
\|N_1(\phi_1)-N_1(\phi_2)\|_{**, D_\e} \le C\gamma\e^\alpha \| \phi_1-\phi_2\|_{L^\infty(D_\e)}
\]
and
\[
\|N_2(\phi_1)-N_2(\phi_2)\|_{*, \partial D_\e} \le C\gamma\e^\alpha \| \phi_1-\phi_2\|_{L^\infty(D_\e)}
\]
where $C$ is independent of $\gamma, \e$.
Thus the operator $A$ is a contraction mapping on $\mathcal{F}_\gamma$
 and therefore a unique fixed point of $A$ exists in this region.

The discussion of the differentiability of $\phi$ with respect to $\xi'$ can be got by the similar proof in \cite{DPM}.
\end{proof}

\section{Variational reduction}\label{s5}
After Lemma \ref{nonlinear},  we already achieved that
\begin{equation}
\begin{cases}
-\Delta (V+\phi)  =  K(\e y)e^{2V+2\phi}, \qquad \qquad  & \text{in }  D_{\e} , \\
\D \F{\partial(V+\phi)}{\partial\bm n} + \e = \kappa(\e y) e^{V+\phi} + c_1\chi_\e Z_{1,\e}, \quad & \text{on }\partial D_{\e}, \medskip\\
\D \int_{D_\e} \chi_\e Z_{1,\e} \phi=0.
\end{cases}
\end{equation}
 To solve our problem, it is sufficient to let $c_1=0$ in (\ref{nonlp}). In this section, we will see $c_1 = 0$ is equivalent to solving a finite dimensional problem.\par

Define the energy functional
\BEN
\widetilde E(U) = \F{1}{2} \int_D |\nabla U|^2 - \F{\e^2}{2} \int_D K e^{2U} + \int_{\partial D}  U -  \e \int_{\partial D}\kappa (x) e^U.
\EEN
Also in terms of $V$,
\BEN
E(V) = \F{1}{2} \int_{D_\e} |\nabla V|^2 - \F{1}{2} \int_{D_e} K(\e y) e^{2V} + \int_{\partial D_\e} \e  V - \int_{\partial D_\e} \kappa (\e y) e^V.
\EEN
Notice that $E(V) = \widetilde E(U) + 2\ln\e\int_{\partial D} \mathrm dx = \widetilde E(U) + 4\pi\ln\e$.

\BP\label{p2}
If $\xi'\in \partial D_\e$ is a tangentially critical point of $F(\xi')=E(V(\xi')+\phi(\xi'))$ on $\partial D_\e$. Then $c_1=0$.
\EP

\begin{proof}
Note that $\partial_{\tau}V(\xi') = Z_{1,\e} + o(1)$ in $L^\infty(D_\e)$, where $\tau$ is the unit tangent vector at $\xi'$.
It is easy to see, since $\|\partial_{\xi'}\phi\|_\infty=o(1)$, that
\begin{align*}
 0=\partial_\tau F(\xi')= E'(V+\phi) (\partial_{\tau}V + \partial_{\tau}\phi) = \int_{D_\e} c_1 \chi_\e Z_{1,\e} (\partial_{\tau}V + \partial_{\tau}\phi)
 =  \int_{D_\e} c_1 \chi_\e Z_{1,\e}^2 + o(1),
\end{align*}
from which the proposition is concluded .
\end{proof}

Next we will compute $E(V+\phi)$, on account of Proposition \ref{p2}.

\BP\label{p1}
It holds that
\BEN
E(V+\phi)=E(V) + O(\e^{2\alpha}).
\EEN
\EP

\begin{proof}
On account of (\ref{nonlp}), we have, for some $t\in(0,1)$, that
\begin{align*}
E(V+\phi) &= E(V) + E'(V)\phi + \F{1}{2}E''(V+t\phi)\phi^2  \\
&= E(V) - \int_{D_\e} R_1 \phi - \int_{\partial D_\e} R_2 \phi \\
&\qquad + \F{1}{2}\int_{D_\e} R_1 \phi + \F{1}{2}\int_{D_\e} N_1(\phi)\phi + \int_{D_\e} K(\e y) e^{2V} (1-e^{2t\phi})\phi^2 \\
&\qquad + \F{1}{2}\int_{\partial D_\e} R_2 \phi + \F{1}{2}\int_{\partial D_\e} N_2(\phi)\phi + \F{1}{2}\int_{\partial D_\e} \kappa(\e y)e^V (1-e^{t\phi})\phi^2.
\end{align*}
From (\ref{error}), it is easily checked that
\begin{align*}
&\left|\int_{D_\e} R_1 \phi \right| \leq C \|R_1\|_{**,D_\e} \|\phi\|_\infty = O(\e^{2\alpha}), && \int_{\partial D_\e} R_2 \phi = O(\e^{2\alpha}),\\
&\int_{D_\e} N_1(\phi)\phi = O(\|\phi\|_\infty^3)=O(\e^{3\alpha}), && \int_{\partial D_\e} N_2(\phi)\phi=O(\e^{3\alpha}).
\end{align*}
Also,
\begin{align*}
\left| \int_{D_\e} K(\e y) e^{2V} (1-e^{2t\phi})\phi^2\right| &\leq C \int_{\mathbb R^2} \F{1}{(\lambda^2 +  |y-\xi'|^2)^2} |\phi|^3=O(\e^{3\alpha}), \\
\left|\int_{\partial D_\e} \kappa(\e y)e^V (1-e^{t\phi})\phi^2\right| &= O(\e^{3\alpha}).
\end{align*}
The proof is complete.
\end{proof}

\section{Energy expansion and the proof of Theorem}\label{s6}
In this section, we will first compute the energy functional  and then prove the main theorem.
\BT
It holds that
\BEN
E(V)= 2\pi\ln\e - 2\pi + 2\pi\ln2  - 2 \pi \ln \left(\kappa(\xi)+\sqrt{K(\xi)+\kappa^2(\xi)}\right) + O(\e^\alpha).
\EEN
\ET


\begin{proof}
From the equation (\ref{3}) of $U$, we may have that
\begin{align*}
  E(V) =&\ \F{1}{2} \int_{D_\e} (-\Delta V)V - \F{1}{2} \int_{D_e} K(\e y) e^{2V} + \F{1}{2} \int_{\partial D_\e}\F{\partial V}{\partial\bm n}V + \int_{\partial D_\e} \e  V - \int_{\partial D_\e} \kappa (\e y) e^V \\
  =&\  \F{1}{2} \int_{D_\e} K(\xi) e^{2(U_0(\e y)+2\ln\e)}V - \F{1}{2} \int_{D_e} K(\e y) e^{2V} + \F{\e}{2} \int_{\partial D_\e}  V - \F{\e d}{2}\int_{\partial D_\e} \F{\lambda^2}{\lambda^2+|y-\xi'|^2}V \\
  &\ + \F{1}{2} \int_{\partial D_\e} \kappa(\xi) e^{U_0(\e y)+2\ln \e}V - \int_{\partial D_\e} \kappa(\e y)  e^V.
\end{align*}
Let us calculate terms by terms. From the computation in the appendix,
we get that
\begin{align*}
 &\ \int_{D_\e} K(\xi) e^{2(U_0(\e y)+2\ln\e)}V \\
 =&\ \int_{D_\e} \F{4\lambda^2}{(\lambda^2+|y-\xi'-\mathfrak D(\xi) \lambda \bm n(\xi)|^2)^2} \left(\ln\F{2\lambda}{\sqrt{K(\xi)}(\lambda^2+|y-\xi'-\mathfrak D(\xi) \lambda \bm n(\xi)|^2)} + H_0(\e y)\right) \\
 =&\ \int_{D_\e} \F{4\lambda^2}{(\lambda^2+|y-\xi'-\mathfrak D(\xi) \lambda \bm n(\xi)|^2)^2} \left( \ln \F{2\lambda}{\sqrt{K(\xi)}} + \ln\F{1}{(\lambda^2+|y-\xi'-\mathfrak D(\xi) \lambda \bm n(\xi)|^2)}  + O(\e^\alpha)\right)\\
 =&\ -2\pi \left[2\sinh^{-1} \mathfrak D(\xi) +1+2\ln\lambda - \F{\mathfrak D(\xi)}{\sqrt{1+\mathfrak D(\xi)^2}} \left(1+2\ln2\lambda +\ln (1+\mathfrak D(\xi)^2)\right)\right]  \\
&\qquad + 2\pi \left(1-\F{\mathfrak D(\xi)}{\sqrt{1+\mathfrak D(\xi)^2}}\right) \ln \F{2\lambda}{\sqrt{K(\xi)}}+ O(\e^\alpha)\\
=&\ -2\pi \left[2\sinh^{-1} \mathfrak D(\xi) - \F{\mathfrak D(\xi)}{\sqrt{1+\mathfrak D(\xi)^2}} \left(2\ln2 +\ln (1+\mathfrak D(\xi)^2)\right)\right]\\
&\qquad + 2\pi \left(1-\F{\mathfrak D(\xi)}{\sqrt{1+\mathfrak D(\xi)^2}}\right) \left( \ln \F{2\lambda}{\sqrt{K(\xi)}}-1 -2\ln\lambda\right)+ O(\e^\alpha).
\end{align*}
Also it holds that
\begin{align*}
  \int_{D_\e} K(\e y) e^{2V} &=  \int_{D_\e} \F{4\lambda^2 K(\e y) e^{2H_0(\e y)}}{K(\xi)(\lambda^2+|y-\xi'-\mathfrak D(\xi) \lambda \bm n(\xi)|^2)^2} \\
  &=\int_{D_\e\cap B_\F{\delta}{\e}(\xi')} \F{ 4\lambda^2}{(\lambda^2+|y-\xi'-\mathfrak D(\xi) \lambda \bm n(\xi)|^2)^2} + O(\e^\alpha)\\
 &= \int_{D_\e} \F{ 4\lambda^2}{(\lambda^2+|y-\xi'-\mathfrak D(\xi) \lambda \bm n(\xi)|^2)^2} + O(\e^\alpha)\\
 &=2\pi \left(1-\F{\mathfrak D(\xi)}{\sqrt{1+\mathfrak D(\xi)^2}}\right)+O(\e^\alpha).
\end{align*}
So we obtain that
\begin{align}
 &\ \F{1}{2}\int_{D_\e} K(\xi) e^{2(U_0(\e y)+2\ln\e)}V - \F{1}{2}\int_{D_\e} K(\e y) e^{2V} \nonumber \\
  =&\ -\pi\left[2\sinh^{-1} \mathfrak D(\xi) - \F{\mathfrak D(\xi)}{\sqrt{1+\mathfrak D(\xi)^2}} \left(2\ln 2+\ln \left(1+\mathfrak D(\xi)^2\right)\right)\right]  \nonumber\\
 &\qquad + \pi  \left (\ln \F{2\lambda}{\sqrt{K(\xi)}}-2 -2 \ln \lambda \right) \left( 1 - \F{\mathfrak D(\xi)}{\sqrt{1+\mathfrak D(\xi)^2}}\right)+ O(\e^\alpha). \label{4}
\end{align}
Again we calculate the boundary terms.
It is checked, due to the appendix, that
\begin{align*}
  &\ \int_{\partial D_\e} \kappa(\xi) e^{U_0(\e y)+2\ln \e}V  \\
  =&\ \int_{\partial  D_\e} \F{2\lambda \mathfrak D(\xi)}{\lambda^2+|y-\xi'-\mathfrak D(\xi) \lambda \bm n(\xi)|^2} \left(\ln\F{2\lambda}{\sqrt{K(\xi)}(\lambda^2+|y-\xi'-\mathfrak D(\xi) \lambda \bm n(\xi)|^2)} + H_0(\e y)\right) \\
  = &\ \int_{\partial D_\e} \F{2\lambda \mathfrak D(\xi)}{\lambda^2+|y-\xi'-\mathfrak D(\xi) \lambda \bm n(\xi)|^2} \left(\ln\F{2\lambda}{\sqrt{K(\xi)}}+\ln\F{1}{(\lambda^2+|y-\xi'-\mathfrak D(\xi) \lambda \bm n(\xi)|^2)} + O(\e^\alpha)\right)\\
  = &\ 2\pi \F{\mathfrak D(\xi)}{\sqrt{1+\mathfrak D(\xi)^2}} \left[\ln\F{2\lambda}{\sqrt{K(\xi)}}-\ln 4 -2\ln\lambda - \ln(1+\mathfrak D(\xi)^2)\right] +O(\e^\alpha).
\end{align*}
and
\begin{multline*}
  \int_{\partial D_\e} \kappa (\e y) e^V = \int_{\partial D_\e\cap B_\F{\delta}{\varepsilon}(\xi')}  \F{2\lambda \kappa(\e y) e^{H_0(\e y)}}{\sqrt{K(\xi)}(\lambda^2+|y-\xi'-\mathfrak D(\xi) \lambda \bm n(\xi)|^2)}  \\
  =  \int_{\partial D_\e} \F{2\lambda \mathfrak D(\xi)}{(\lambda^2+|y-\xi'-\mathfrak D(\xi) \lambda \bm n(\xi)|^2)} + O(\e^\alpha)  =2\pi  \F{\mathfrak D(\xi)}{\sqrt{1+\mathfrak D(\xi)^2}}  + O(\e^\alpha).
\end{multline*}
It then holds that
\begin{multline} \label{5}
   \F{1}{2} \int_{\partial D_\e} \kappa(\xi) e^{U_0(\e y)+2\ln \e}V - \int_{\partial D_\e} \kappa (\e y) e^V  \\
  =\pi  \F{\mathfrak D(\xi)}{\sqrt{1+\mathfrak D(\xi)^2}} \left[ \ln \F{2\lambda}{\sqrt{K(\xi)}} - 2\ln 2 -2-2\ln\lambda - \ln(1+\mathfrak D(\xi)^2)\right]+ O(\e^\alpha).
\end{multline}
Obviously we see that
\BE\label{6}
-\F{\e d}{2}\int_{\partial D_\e} \F{\lambda^2}{\lambda^2+|y-\xi'|^2}V= O(\e).
\EE
The last term $ \F{\e}{2} \int_{\partial D_\e}  V $ is a nonlocal term, which is checked that
\begin{align}
 \F{\e}{2} \int_{\partial D_\e}  V =&\ \F{\e}{2} \int_{\partial D_\e}  \left(\ln\F{2\lambda}{\sqrt{K(\xi)}(\lambda^2+|y-\xi'-\mathfrak D(\xi) \lambda \bm n(\xi)|^2)} + H_0(\e y)\right) \nonumber \\
 =&\ \pi\ln \F{2\lambda}{\sqrt{K(\xi)}}- \F{\e}{2} \int_{\partial D_\e} \ln(\lambda^2+|y-\xi'-\mathfrak D(\xi) \lambda \bm n(\xi)|^2)  + O(\e^\alpha)  \nonumber \\
 =&\ \pi\ln \F{2\lambda}{\sqrt{K(\xi)}} + 2\pi\ln\e  +  O(\e^\alpha). \label{7}
\end{align}
From (\ref{4})--(\ref{7}), we achieve that
\begin{multline*}
 E(V) = 2\pi\ln\e - 2\pi  +2\pi \left( \ln\F{2\lambda}{\sqrt{K(\xi)}} -\sinh^{-1}\mathfrak D(\xi) - \ln\lambda \right) +  O(\e^\alpha) \\
 = 2\pi\ln\e - 2\pi  +2\pi\ln 2 - 2\pi \ln \sqrt{K(\xi)} - 2\pi \sinh^{-1}\mathfrak D(\xi) +O(\e^\alpha)\\
 =2\pi\ln\e - 2\pi + 2\pi\ln2 - 2\pi \ln \left( \kappa(\xi) + \sqrt{K(\xi)+\kappa(\xi)^2}\right) +O(\e^\alpha).
\end{multline*}
The proof is complete.
\end{proof}

%

\begin{proof}[Proof of Theorem \ref{T1}]
Since $\xi_*$ is a local extremum  point of $\kappa(\xi) + \sqrt{K(\xi)+\kappa(\xi)^2}$ on the boundary,  there is an open neighborhood $\Gamma\subset\partial D$ of $\xi_*$ such that
\BEN
\kappa(\xi) + \sqrt{K(\xi)+\kappa(\xi)^2} <\ (\text{or } >)\ \kappa(\xi_*) + \sqrt{K(\xi_*)+\kappa(\xi_*)^2}, \qquad \forall~\xi\in\partial\Gamma.
\EEN
Recall again that
\BEN
E(V+\phi)[\xi] = 2\pi\ln\e - 2\pi + 2\pi\ln2 - 2\pi \ln \left( \kappa(\xi) + \sqrt{K(\xi)+\kappa(\xi)^2}\right) +O(\e^\alpha).
\EEN
So we have, from the continuity, that
\BEN
E(V+\phi)[\xi] < \ (\text{or } >)\  E(V+\phi)[\xi_*], \qquad \forall~\xi\in\partial\Gamma.
\EEN
Therefore $\xi_*\in\partial D$ must be a tangentially critical point of $E(V+\phi)[\xi]$. The existence result is then from Proposition \ref{p2}. The remaining identity is just due to the Gauss-Bonnet theorem.
\end{proof}

%
%

\section{Appendix}\label{s7}
Some useful computation is given in this section.

\BL
It holds that, for $\lambda>0$,
\BEN
\int_{D_\e} \F{4\lambda^2}{(\lambda^2+|y-\xi'-\mathfrak D(\xi) \lambda \bm n(\xi)|^2)^2}\mathrm dy = 2\pi \left(1-\F{\mathfrak D(\xi)}{\sqrt{1+\mathfrak D(\xi)^2}}\right) + O(\e).
\EEN
\EL

\begin{proof}
Without loss of generality, $\xi$ is assumed to be the origin owing to the symmetry. For a small and fixed $\delta>0$, we have that
\begin{align*}
 &\ \int_{D_\e} \F{4\lambda^2}{(\lambda^2+|y-\mathfrak D(0) \lambda \bm n(0)|^2)^2}\mathrm dy
=  \int_{D_{\lambda \e}} \F{4}{(1+|z-\mathfrak D(0)\bm n(0)|^2)^2}\mathrm dz  \\
=&\ \int_{D_{\lambda \e} \cap B^+_{\delta/\lambda\e}(0)} \F{4}{(1+|z-\mathfrak D(0)\bm n(0)|^2)^2}\mathrm dz + O(\e^2) \\
=&\ \int_{ B^+_{\delta/\lambda\e}(0)} \F{4}{(1+|z-\mathfrak D(0)\bm n(0)|^2)^2}\mathrm dz - \int_{D_{\lambda \e}^c \cap B^+_{\delta/\lambda\e}(0)} \F{4}{(1+|z-\mathfrak D(0)\bm n(0)|^2)^2}\mathrm dz + O(\e^2) \\
=&\ \int_{\mathbb R^2_+} \F{4}{(1+|z-\mathfrak D(0)\bm n(0)|^2)^2}\mathrm dz + O(\e) = 2\pi \left(1-\F{\mathfrak D(0)}{\sqrt{1+\mathfrak D(0)^2}}\right) + O(\e) ,
\end{align*}
where the last but one equality is due to
\begin{multline*}
  \int_{D_{\lambda \e}^c \cap B^+_{\delta/\lambda\e}(0)} \F{4}{(1+|z-\mathfrak D(0)\bm n(0)|^2)^2}\mathrm dz \leq 2\int_0^\F{\delta}{\lambda\e} \mathrm d z_1 \int_0^{\F{1}{\lambda\e}-\sqrt{\F{1}{\lambda^2\e^2}-z_1^2}}\F{4}{(1+z_1^2+z_2^2)^2}\mathrm dz_2 \\
  \leq 2\int_0^\F{\delta}{\lambda\e} \F{4}{(1+z_1^2)^2}\left(\F{1}{\lambda\e}-\sqrt{\F{1}{\lambda^2\e^2}-z_1^2}\right) \mathrm d z_1 = O(\lambda \e). \qedhere
\end{multline*}
\end{proof}

Similarly, the following calculation holds too.
\BL
We have
\begin{multline*}
\int_{D_\e} \F{4\lambda^2}{(\lambda^2+|y-\xi'-\mathfrak D(\xi) \lambda \bm n(\xi)|^2)^2}\ln\F{1}{(\lambda^2+|y-\xi'-\mathfrak D(\xi) \lambda \bm n(\xi)|^2)}\mathrm dy \\
= -2\pi \left[2\sinh^{-1} \mathfrak D(\xi) +1+2\ln\lambda - \F{\mathfrak D(\xi)}{\sqrt{1+\mathfrak D(\xi)^2}} \left(1+2\ln2\lambda +\ln (1+\mathfrak D(\xi)^2)\right)\right] +O(\e).
\end{multline*}
\EL

\begin{proof}
It holds that
\begin{align*}
 &\ \int_{D_\e} \F{4\lambda^2}{(\lambda^2+|y-\mathfrak D(0) \lambda \bm n(0)|^2)^2}\ln\F{1}{(\lambda^2+|y-\mathfrak D(0) \lambda \bm n(0)|^2)}\mathrm dy \\
=&\ \int_{D_{\lambda \e}} \F{4}{(1+|z-\mathfrak D(0)\bm n(0)|^2)^2} \left[\ln \F{1}{(1+|z-\mathfrak D(0)\bm n(0)|^2)} -2\ln\lambda \right]\mathrm dz \\
=&\ \int_{\mathbb R^2_+} \F{4}{(1+|z-\mathfrak D(0)\bm n(0)|^2)^2}\ln \F{1}{(1+|z-\mathfrak D(0)\bm n(0)|^2)}\mathrm dz - 2\ln\lambda \int_{\mathbb R^2_+} \F{4}{(1+|z-\mathfrak D(0)\bm n(0)|^2)^2}\mathrm dz\\
&\qquad\qquad  + O(\e)\\
=&\ -2\pi \left[2\sinh^{-1} \mathfrak D(0) +1- \F{\mathfrak D(0)}{\sqrt{1+\mathfrak D(0)^2}} \left(1+\ln 4+\ln (1+\mathfrak D(0)^2)\right)\right]\\
&\qquad\qquad -4\pi \ln \lambda \left(1-\F{\mathfrak D(0)}{\sqrt{1+\mathfrak D(0)^2}}\right) + O(\e) \\
=&\ -2\pi \left[2\sinh^{-1} \mathfrak D(0) +1+2\ln\lambda - \F{\mathfrak D(0)}{\sqrt{1+\mathfrak D(0)^2}} \left(1+2\ln2\lambda +\ln (1+\mathfrak D(0)^2)\right)\right] +O(\e). \qedhere
\end{align*}
\end{proof}

The next are about the boundary terms.
\BL
We have
\BEN
\int_{\partial D_\e} \F{2\lambda \mathfrak D(\xi)}{\lambda^2+|y-\xi'-\mathfrak D(\xi) \lambda \bm n(\xi)|^2}
 = 2\pi \F{\mathfrak D(\xi)}{\sqrt{1+\mathfrak D(\xi)^2}} +O(\e).
\EEN
\EL
\begin{proof}
Calculation shows that
\begin{align*}
 &\ \int_{\partial D_\e} \F{2\lambda \mathfrak D(0)}{\lambda^2+|y-\mathfrak D(0) \lambda \bm n(0)|^2} \mathrm dy
 = \int_{\partial D_{\lambda\e}\cap B_{\F{\delta}{\lambda\e}}} \F{2 \mathfrak D(0)}{1+|z-\mathfrak D(0) \bm n(0)|^2} \mathrm dz + O(\e) \\
 =&\ \int_{-\F{\tilde\delta}{\lambda\e}}^{\F{\tilde\delta}{\lambda\e}} \F{2 \mathfrak D(0)}{1+z_1^2+\left(\F{1}{\lambda\e}-\F{1}{\lambda\e}\sqrt{1-\lambda^2\e^2z_1^2}+\mathfrak D(0)\right)^2} \F{1}{\sqrt{1-\lambda^2\e^2z_1^2}} \mathrm dz_1 + O(\e) \\
 =&\ \int_{-\F{\tilde\delta}{\lambda\e}}^{\F{\tilde\delta}{\lambda\e}} \F{2\mathfrak D(0)}{1+z_1^2+\mathfrak D(0)^2} \mathrm dz_1 + O(\e)  = \int_{\partial\mathbb R^2_+}\F{2 \mathfrak D(0)}{1+|z-\mathfrak D(0) \bm n(0)|^2}  + O(\e)\\
 =&\   2\pi \F{\mathfrak D(0)}{\sqrt{1+\mathfrak D(0)^2}}+ O(\e). \qedhere
\end{align*}
\end{proof}

\BL
It holds that
\begin{multline*}
  \int_{\partial D_\e} \F{2\lambda \mathfrak D(\xi)}{\lambda^2+|y-\xi'-\mathfrak D(\xi) \lambda \bm n(\xi)|^2} \ln\F{1}{(\lambda^2+|y-\xi'-\mathfrak D(\xi) \lambda \bm n(\xi)|^2)}\mathrm dy \\
  =- 2\pi \F{\mathfrak D(\xi)}{\sqrt{1+\mathfrak D(\xi)^2}}\left[ \ln 4 + \ln (1+\mathfrak D(\xi)^2) + 2\ln\lambda \right]+ O(\e^{\alpha}).
\end{multline*}
\EL

\begin{proof}
Direct computation shows that
\begin{align*}
 &\ \int_{\partial D_\e} \F{2\lambda \mathfrak D(\xi)}{\lambda^2+|y-\xi'-\mathfrak D(\xi) \lambda \bm n(\xi)|^2} \ln\F{1}{(\lambda^2+|y-\xi'-\mathfrak D(\xi) \lambda \bm n(\xi)|^2)}\mathrm dy \\
 =&\ \int_{\partial D_{\lambda\e}\cap B_{\F{\delta}{\lambda\e}}} \F{2 \mathfrak D(0)}{1+|z-\mathfrak D(0) \bm n(0)|^2} \left[\ln\F{1}{(1+|z-\mathfrak D(0)  \bm n(0)|^2)} - 2\ln\lambda \right]\mathrm dz + O(\e^{\alpha}) \\
 =&\ \int_{-\F{\tilde\delta}{\lambda\e}}^{\F{\tilde\delta}{\lambda\e}} \F{2 \mathfrak D(0)}{1+z_1^2+\left(\F{1}{\lambda\e}-\F{1}{\lambda\e}\sqrt{1-\lambda^2\e^2z_1^2}+\mathfrak D(0)\right)^2}\left[\ln\F{1}{1+z_1^2+\left(\F{1}{\lambda\e}-\F{1}{\lambda\e}\sqrt{1-\lambda^2\e^2z_1^2}+\mathfrak D(0)\right)^2} - 2\ln\lambda \right] \\
 &\ \qquad\qquad \cdot \F{1}{\sqrt{1-\lambda^2\e^2z_1^2}} \mathrm dz_1 + O(\e^{\alpha}) \\
 =&\ \int_{-\F{\tilde\delta}{\lambda\e}}^{\F{\tilde\delta}{\lambda\e}} \F{2\mathfrak D(0)}{1+z_1^2+\mathfrak D(0)^2}\left[\ln\F{1}{1+z_1^2+\mathfrak D(0)^2}-2\ln\lambda\right] \mathrm dz_1 + O(\e^{\alpha}) \\
 =&\ \int_{\partial\mathbb R^2_+}\F{2 \mathfrak D(0)}{1+|z-\mathfrak D(0) \bm n(0)|^2} \left[\ln\F{1}{(1+|z-\mathfrak D(0)  \bm n(0)|^2)} - 2\ln\lambda \right]\mathrm dz + O(\e^{\alpha})\\
 =&\  - 2\pi \F{\mathfrak D(\xi)}{\sqrt{1+\mathfrak D(\xi)^2}}\left[ \ln 4 + \ln (1+\mathfrak D(\xi)^2) + 2\ln\lambda \right]+ O(\e^{\alpha}).\qedhere
\end{align*}
\end{proof}

\BL
It holds that
\BEN
\e \int_{\partial D_\e} \ln(\lambda^2+|y-\xi'-\mathfrak D(\xi) \lambda \bm n(\xi)|^2) = -4\pi\ln\e +O(\e).
\EEN
\EL
\begin{proof}
We may check that
\begin{align*}
&\ \e \int_{\partial D_\e} \ln(\lambda^2+|y-\mathfrak D(0) \lambda \bm n(0)|^2)\mathrm dy = \lambda\e\int_{\partial D_{\lambda\e}}\left[\ln\lambda^2 + \ln(1+|z-\mathfrak D(0)\bm n(0)|^2)\right] \mathrm d z\\
=&\ 4\pi\ln\lambda +  \lambda\e\int_{\partial D_{\lambda\e}}\ln\left(1+|z-\mathfrak D(0)\bm n(0)|^2\right)\mathrm dz \\
=&\ 4\pi\ln\lambda +2\pi  \left(-2 \ln (\lambda\e)+ O(\e)\right) = -4\pi \ln\e + O(\e). \qedhere
\end{align*}
\end{proof}

\end{document}